\RequirePackage{fix-cm}
\documentclass[smallextended]{svjour3}       
\usepackage[utf8]{inputenc}

\usepackage{tabularx,ragged2e,booktabs,caption}
\newcolumntype{C}[1]{>{\Centering}m{#1}}

\smartqed  
\usepackage{xspace}
\usepackage{amsfonts}
\usepackage{amsmath}
\usepackage{enumerate}
\usepackage{framed}
\usepackage{graphicx}
\usepackage{lscape}
\usepackage{bm}
\usepackage{color}
\usepackage{multicol}
\usepackage{cite}
\usepackage[utf8]{inputenc}

\usepackage{tabularx,ragged2e,booktabs,caption}
\newcolumntype{C}[1]{>{\Centering}m{#1}}

\usepackage{amsmath}
\usepackage{amsfonts}
\usepackage{amssymb}
\usepackage{tikz}
\usepackage{subfig}
\usepackage{float}
\usepackage{slashbox}
\usepackage{graphicx}
\usepackage[english]{babel}
\usepackage{geometry}                
\usepackage[parfill]{parskip}    
\usepackage{epstopdf}
\usepackage{psfrag}

\usepackage[sort]{natbib}
\renewcommand{\eqref}[1]{{(\ref{#1})}}

\newcommand{\N}{\mathbb{N}}
\newcommand{\M}{\mathbf{M}}
\newcommand{\C}{\mathcal{C}}
\newcommand{\V}{\mathcal{V}}
\newcommand{\R}{\mathcal{R}}
\newcommand{\F}{\mathcal{F}}
\newcommand{\E}{\mathcal{E}}
\newcommand{\W}{\mathcal{W}}

\usepackage{graphicx}
\begin{document}

\title{A fitted L-Multi-point Flux Approximation method for pricing options}
%


\titlerunning{A fitted L-Multi-point Flux Approximation method for pricing options}        

\author{ Rock Stephane Koffi,
Antoine Tambue
}

\institute{A. Tambue (Corresponding author) \at
	Western Norway University of Applied Sciences,  Inndalsveien 28, 5063  Bergen, Norway,\\
	The African Institute for Mathematical Sciences(AIMS), 6-8 Melrose Road, Muizenberg 7945, South Africa\\
	Center for Research in Computational and Applied Mechanics (CERECAM), and Department of Mathematics and Applied Mathematics, University of Cape Town, 7701 Rondebosch, South Africa.\\
	Tel.: +47 55 58 70 06, \email{antonio@aims.ac.za, antoine.tambue@hvl.no, tambuea@gmail.com}\\
   \and
   R.S. Koffi \at
            The African Institute for Mathematical Sciences(AIMS) of South Africa\\
            Department of Mathematics and Applied Mathematics, University of Cape Town, 7701 Rondebosch, South Africa\\
            \email{rock@aims.ac.za}          
}
\date{Received: date / Accepted: date}

\maketitle

\begin{abstract}
In this paper, we introduce a special kind of finite volume method called Multi-Point Flux Approximation method (MPFA) to price European and American options in two dimensional domain. We focus on the L-MPFA method for space discretization of the diffusion term of Black-Scholes operator. The degeneracy of the Black Scholes operator is tackled using the fitted finite volume method. This combination of fitted finite volume method and L-MPFA method coupled to upwind methods  gives us a novel scheme called the fitted L-MPFA method. Numerical experiments  show the accuracy of the  novel fitted L-MPFA method comparing to well known schemes for pricing options.
\end{abstract}

	
\section{Introduction}
\label{intro}
In finance, an option is a contract which gives to the holder the right but not the obligation to buy (call) or to sell (put) an asset at a a specific price (strike) at a certain date in the future (expiry date). We have two main types of options which are European and American options. European options are options that can be exercised only at expiry date while American options can be exercised anytime before the expiry date. This flexibility of exercising American options leads to solve an optimal stopping time problem in the Black-Scholes framework which incorporate the early exercise.
Many studies focused on the pricing problem of American options were conducted and the linear complementary problem approach was  quite popular for pricing American options (see \cite{kovalov2007pricing,zhang2009power,wang2006power,topper2005financial}). This approach brings us to solve linear complementary problem stated as follows (see \cite{topper2005financial}):

\begin{align}
\left\lbrace
\begin{array}{lcr}
\mathcal{L}V & \geq & 0\\
& & \\
V-V^{\star} & \geq & 0\\
&    & \\
\mathcal{L}V \cdot \big(V-V^{\star}\big) &  = & 0 
\end{array}
\right.
\end{align}

%

where $\mathcal{L}$ is the following Black-Scholes operator
\begin{eqnarray}
\label{BS-op}
\mathcal{L}V=\dfrac{\partial V}{\partial \tau }-\frac{1}{2}\sum_{i,j=1}^n \sigma_i\sigma_j\rho_{ij}S_iS_j\frac{\partial^2V}{\partial S_i \partial S_j}-r\sum_{i=1}^nS_i\frac{\partial V}{\partial S_i}+rV
\end{eqnarray}

with $r$ is the risk free interest, $V$ is the option  value at time $\tau$, $V^{\star}$ is the payoff, $\tau=T-t$ with $t$ and $T$ respectively the instantaneous and  maturity time. 
For $i,j=1,\ldots,n,\,S_i~$ represents the asset $i$ price, $\sigma_i$  represents the volatility of asset $i$,  $\rho_{ij}$  represents the correlation between the assets $i$ and $j$.\\
Furthermore, \cite{wang2006power} proposed a power penalty method to solve the linear complementary problem for pricing  American options. The power penalty problem is formulated as follows:
\begin{equation}
\label{penpow}
\mathcal{L}V+\beta\big[V^{\star}-V\big]_{+}^{1/k}=0
\end{equation}
where $\beta$ is penalty parameter and $k$ is the power of the method.
Let us notice that, when we take the  penalty parameter $\beta=0$ in (\ref{penpow}), we get the Black-scholes Equation for pricing European options, with the operator $\mathcal{L}$ defined in (\ref{BS-op}). However, the power penalty problem (\ref{penpow}) can not be solved analytically, therefore numerical methods are required for its resolution.
Nevertheless, the Black-Scholes operator (\ref{BS-op}) is degenerated when the stock price approaches zero. This degeneracy can affect the accuracy of the numerical method used for the resolution. To tackle this problem, several methods have been proposed. The fitted finite volume method, proposed by S.Wang in \cite{wang2004novel} whereby a rigorous proof of convergence is provided, appears to be more attractive. Moreover, the fitted finite volume method has been used for the resolution of the two dimensional second order Black Scholes PDE followed by the convergence proof in \cite{huang2006fitted}. In spite of the fact that the fitted finite volume methods perform well for the resolution of the Black Scholes PDE, they are only of order 1 with respect to asset price variable. Besides, the fitted O-Multi-Point Flux Approximation (O-MPFA) method has been proposed in \cite{rock2019mpfa} to overcome the degeneracy problem of the Black-Scholes PDE. It has been shown that the O-MPFA is more accurate than the classical fitted finite volume method by \cite{wang2004novel}. However,  the O-MPFA  is heavy (9 points stencil method) and for more general grids, the convergence rate of the O-MPFA method  may reduce (see \cite{aavatsmark2007multipoint}).
%
%

 In this  paper, we focus on the L-MPFA method which is based on the approximation of a linear function gradient defined over a given triangle and the continuity of flux through the edges of this triangle.\\
Indeed, the L-MPFA method is a 7 points stencil method while the O-MPFA is a 9 points stencil method.  This shows that the O-MPFA method can be computationally more expensive than the L-MPFA method. 
Moreover, for more general grids, the order reduction in convergence rate is larger for the O-MPFA than the L-MPFA (see \cite{aavatsmark2002introduction}).
Thereby, to approximate the solution of the second order Black Scholes operator, we couple the L-MPFA method with the upwind methods (first and second order).
Besides, the degeneracy of the Black Scholes operator (\ref{BS-op}) is handled by the fitted finite volume ,\cite{wang2004novel}, when the stock price is approaching zero. The L-MPFA method coupled with the upwind methods ($1^{st}$ and $2^{nd}$ order) is used to approximate the solution of (\ref{penpow}) when the Black operator is not degenerated. We call fitted L-MPFA method that combination of the L-MPFA method and the fitted finite volume method.
Numerical simulations show that  the new fitted L-MPFA method is more accurate than the  fitted O-MPFA  method developed in \cite{rock2019mpfa}  and the standard fitted finite volume method developed in \cite{huang2006fitted}.

The paper is structured as follows. In  section 2, we present the power penalty problem with the corresponding  initial and boundary conditions.
  The spatial discretization of the linear operator is developed in the section 3.  Details on the L-MPFA method of the diffusion term discretization are provided. The convection term is discretized using the upwind methods ($1^{st}$ and $2^{nd}$ method).
   At the end of the section 3,  the novel fitted MPFA method is provided. The $\theta-$ Euler method is used for the time discretization method in the section 4.
    Numerical experiments are presented for the different numerical methods are presented in the section 5. 
    The conclusions of our study are drawn in the last section.

\section{Formulation of the problem}

\subsection{Option with two underlying assets}

Pricing an American option with 2 underlying assets leads to solve the following power penalty problem:

\begin{equation}
\label{penpow1}
\mathcal{L}V+\beta\big[V^{\star}-V\big]_{+}^{1/k}=0
\end{equation}

where the Black-Scholes operator is defined as:

\begin{eqnarray}
\label{BS-op1}
\mathcal{L}V=\frac{\partial V}{\partial \tau}-\frac{1}{2} \sigma^2_1 x^2 \frac{\partial^2 V}{\partial x^2} -\rho \sigma_1 \sigma_2 xy\frac{\partial^2 V}{\partial x \partial y} - \frac{1}{2} \sigma_2^2 y^2 \frac{\partial^2 V}{\partial y^2}-rx\frac{\partial V}{\partial x} -ry\frac{\partial V}{\partial y}+rV
\end{eqnarray}

with the following initial and boundary conditions

\begin{align}
\left\lbrace
\begin{array}{l}
V(x,y,0) =  V^*(x,y)  =  \max\left(K-\alpha_1 x-\alpha_2 y,0\right) \\
 \\
V(0,y,\tau,)=V(x,0,\tau,)  =  K  \\
 \\
\lim_{x,y\longrightarrow + \infty}  V(x,y,\tau)  =   0  
\end{array}
\right.
\end{align}

with $K$ the strike price, $V^*$ the payoff for basket options and $\alpha_i,~i=1,2$, are weights such that $\alpha_1+\alpha_2=1$.\\
\\
When the penalty parameter $\beta=0$ in (\ref{penpow1}), we get the Black-Scholes Partial Differential Equation for pricing European options with the corresponding initial and boundary conditions

\begin{align}
\left\lbrace
\begin{array}{l}
V(x,y,0)=V^*(x,y)=\max\left(K-\alpha_1 x-\alpha_2 y,0\right) \\
\\
V(0,y,\tau)=V(x,0,\tau)=0\\
\\
\lim_{x\longrightarrow\infty}V(x,y,\tau)=x-Ke^{-r\tau}\\ 
\\
\lim_{y\longrightarrow\infty}V(x,y,\tau)=y-Ke^{-r\tau} 
\end{array}\right.
\end{align}

In order to apply the finite volume method, it is convenient to re-write the Black-Scholes operator \eqref{BS-op1} in the following divergence form

\begin{equation}
\label{conservation}
\frac{\partial V}{\partial \tau }-\nabla \cdot(\M\nabla V)-\nabla (f V)-\lambda V+\beta\big[V^{\star}-V\big]_{+}^{1/k}=0
\end{equation}

where

\begin{eqnarray*}
& & \M=\frac{1}{2}\left(\begin{array}{lr}
\sigma_1^2 x^2 & \rho\sigma_1\sigma_2xy \\
   & \\
\rho\sigma_1\sigma_2xy & \sigma_2^2 y^2
\end{array}\right),
f=\left(\begin{array}{c}
(r-\sigma_1^2-\frac{1}{2}\rho\sigma_1\sigma_2)x \\ \\ (r-\sigma_2^2-\frac{1}{2}\rho\sigma_1\sigma_2)y
\end{array}\right)\\
& & \\
&  &\\
&   &~~~~~~~~~~~~~~~~~~~~~~~\lambda = -3r+\sigma_1^2+\sigma^2_2+\rho\sigma_1\sigma_2
\end{eqnarray*}

\subsection{Finite volume method}

Let us consider the new domain $\Omega$ of study by truncating $D$ such that  $\Omega=I_x\times I_y\times [0,T]$ where $I_x=[0,x_{\max}]$ and $I_y=[0,y_{\max}]$.
In  the sequel  of  this work, the Black-Scholes operator \eqref{BS-op1} is considered  over the truncated domain $\Omega$. 

At $x=x_{\max}$ and $y=y_{\max}$, the linear boundary condition will be applied \cite{huang2006fitted}.
The intervals $I_x$ and $I_y$ will be subdivided into  $N$ in the following way \cite{huang2006fitted} \cite{huang2009convergence} without  loss the generality.

\begin{eqnarray}
I_{x_i}=[x_{i-1};x_i], \,\,I_{y_j}=[y_{j-1};y_j]\quad i,j =1,...,N+1.
\end{eqnarray}

Let us set the mid-points $x_{i-\frac{1}{2}}$ and $y_{j-\frac{1}{2}}$ as follows
\begin{eqnarray}
x_{i-\frac{1}{2}}=\frac{x_{i-1}+x_i}{2} \,\,\, y_{j-\frac{1}{2}}=\frac{y_{j-1}+y_j}{2} \qquad i,j =1,...,N,
\end{eqnarray}
with  $h_i=x_{i+\frac{1}{2}}-x_{i-\frac{1}{2}},~~~l_j=y_{j+\frac{1}{2}}-y_{j-\frac{1}{2}}$~~~~~and

\begin{eqnarray*}
x_{-\frac{1}{2}}=x_0=0, \qquad  x_{N+\frac{3}{2}}=x_{N+1}=x_{\max},\,  y_{-\frac{1}{2}}=y_0=0 \,, \;\,y_{N+\frac{3}{2}}=y_{N+1}=y_{\max}.
\end{eqnarray*}

For $i,j=1,\ldots,N+1$,  we denote by $\C_{ij}=[x_{i-\frac{1}{2}};x_{i+\frac{1}{2}}]\times[y_{j-\frac{1}{2}};y_{j+\frac{1}{2}}]$  a  control volume associated our subdivision.

\begin{figure}[!h]
\centering
\begin{tikzpicture}[scale=0.4]
\draw[black,thick] (0,0)--(0,16)--(16,16)--(16,0)--(0,0);
\draw[black,thick] (4,0)--(4,16);
\draw[black,thick] (8,0)--(8,16);
\draw[black,thick] (12,0)--(12,16);
\draw[black,thick] (0,4)--(16,4);
\draw[black,thick] (0,8)--(16,8);
\draw[black,thick] (0,12)--(16,12);
\draw[red,thick] (4,4)--(4,12)--(12,12)--(12,4)--(4,4);
\draw[red,thick] (8,4)--(8,12);
\draw[red,thick] (4,8)--(12,8);
\draw[blue,thick] (6,6)--(6,10)--(10,10)--(10,6)--(6,6);
\draw[red,fill=red] (4,4) circle (0.1);
\draw[red,fill=red] (4,8) circle (0.1);
\draw[red,fill=red] (4,12) circle (0.1);
\draw[red,fill=red] (8,4) circle (0.1);
\draw[red,fill=red] (8,8) circle (0.1);
\draw[red,fill=red] (8,12) circle (0.1);
\draw[red,fill=red] (12,4) circle (0.1);
\draw[red,fill=red] (12,8) circle (0.1);
\draw[red,fill=red] (12,12) circle (0.1);
\node[below,red] at (8.7,8){$(x_i,y_j)$};
\node[above,blue] at (6.7,6){$\C_{ij}$};
\end{tikzpicture}
\caption{Control volume}
\end{figure}
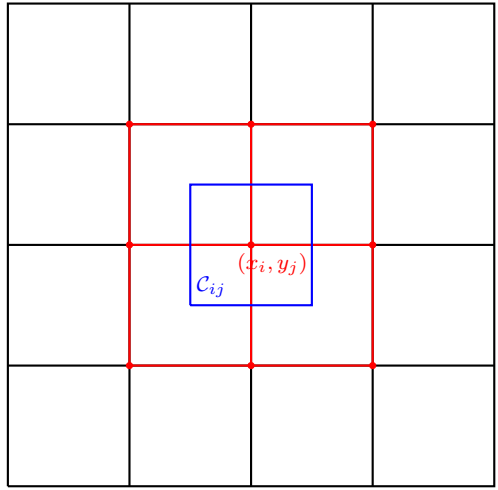
Note  that for $i,j=1,\ldots,N$ the control volume $\C_{ij}$ is the area surrounding the grid point $(x_i,y_j)$.

Our goal  is to approximate  the option  function $V$ at $(x_i,y_j)$  \footnote{center of the control volume $\C_{ij}$}  by a function denoted $\V$.

 The matrix $\M$ in \eqref{conservation} will be replaced by its average value  in each control  volume

\begin{equation}
\M^{ij}=\frac{1}{\mathrm {meas}(\C_{ij})}\int_{\C_{ij}}\M dxdy,\,\,\,  i,j=1,...,N.
\end{equation}
where $~\mathrm{meas}(\C_{ij})$ is the measure of $\C_{ij}$.
Thereby, we have  
\begin{equation*}
M^{ij}= \left[\begin{array}{lcr}
\frac{\sigma_1^2}{6}\frac{x_{i+\frac{1}{2}}^3-x_{i-\frac{1}{2}}^3}{x_{i+\frac{1}{2}}-x_{i-\frac{1}{2}}}&    &\frac{\rho\sigma_1\sigma_2}{8}(x_{i+\frac{1}{2}}+x_{i-\frac{1}{2}})(y_{j+\frac{1}{2}}+y_{j-\frac{1}{2}}) \\
                   &        &  \\
\frac{\rho\sigma_1\sigma_2}{8}(x_{i+\frac{1}{2}}+x_{i-\frac{1}{2}})(y_{j+\frac{1}{2}}+y_{j-\frac{1}{2}})  &    &  \frac{\sigma_2^2}{6}\frac{y_{j+\frac{1}{2}}^3 -y_{j-\frac{1}{2}}^3}{y_{j+\frac{1}{2}}-y_{j-\frac{1}{2}}}
\end{array}\right]
\end{equation*}

Now let us consider the divergence form given in \eqref{conservation}. According to the finite volume method, we integrate the partial differential equation \eqref{conservation} over each control volume $\C_{ij}$. For $i,j=1,...,N~$, we have:

\begin{align}
\label{eqfinvol}
\int_{\C_{ij}}\frac{\partial{\V}}{\partial{\tau}}d\C-\int_{\C_{ij}}\nabla \cdot(\M\nabla \V) d\C-\int_{\C_{ij}}\nabla (f \V)d\C
-\int_{\C_{ij}}\lambda \V d\C+\int_{\C_{ij}}\beta\big[V^{\star}-\V\big]_{+}^{1/k}d\C=0
\end{align}

The next section will be dedicated to spatial discretization of equation \eqref{eqfinvol}. For the term in the  left hand side of the equality sign and for the last one in the right hand side of \eqref{eqfinvol}, we use the mid-point quadrature rule for their approximation as follows: 
\begin{eqnarray}
\int_{{\mathcal{C}}_{ij}}\frac{\partial{\V}}{\partial{\tau}}d\C \approx  meas(\C_{ij})\frac{\partial{\V}}{\partial{\tau}}(x_i,y_j,\tau) \approx  h_il_j\frac{\partial{\V_{ij}}}{\partial{\tau}}
\end{eqnarray}

\begin{eqnarray}
\label{linearterm}
\int_{{\mathcal{C}}_{ij}}\lambda \V d\C \approx meas(\C_{ij})\lambda \V(x_i,y_j,\tau) \approx  h_il_j\lambda \V_{ij}
\end{eqnarray}

\begin{eqnarray}
\label{powpenterm}
\int_{{\mathcal{C}}_{ij}}\beta\big[V^{\star}-\V\big]_{+}^{1/k} d\C \approx meas(\C_{ij})\beta\big[V^{\star}-\V\big]_{+}^{1/k} \approx  h_il_j\beta\big[V_i^{\star}-\V_{ij}\big]_{+}^{1/k}
\end{eqnarray}

the convection term 

\begin{equation}
\label{convectionterm}
\int_{\C_{ij}}\nabla (f \V)d\C
\end{equation}

of \eqref{eqfinvol} will be approximated using the upwind methods (first or second order).

The diffusion term 

\begin{equation}
\label{diffusionterm}
\int_{\C_{ij}}\nabla \cdot(\M\nabla \V) d\C
\end{equation}
of \eqref{eqfinvol} will be approximated using the \textbf{Multi-Point Flux Approximation} (MPFA) \textbf{L}-method or the \textbf{fitted multi-point flux approximation L-} method. More details about these methods will be given in the next section. \\

\section{Space discretization}

The spatial discretization of \eqref{conservation} consists in approximating all terms in \eqref{eqfinvol} over the control volumes of the study domain. 

\subsection{Discretization of the diffusion term}

Let us start by applying the divergence theorem to the diffusion term \eqref{diffusionterm} as follows, for$~~i,j=1,...,N$ :

\begin{equation}
\label{diffterm-disc}
\F^{ij}=\int_{\C_{ij}}\nabla \cdot(\M^{ij}\nabla \V)=\int_{\partial \C_{ij}}(\M^{ij}\nabla \V)\cdot\vec{n}d\partial\C
\end{equation}

where $\vec{n}$ is the outward vector from the control volume.\\
\\
Now, we can apply the so-called L-$\textbf{Multi-Point Flux Aprroximation(MPFA)}$ method   to approximate the integral defined in \eqref{diffterm-disc}.

\subsubsection{L -Multi-Point Flux Approximation (L-MPFA) method }
The L-MPFA method takes its name from the fact that the curve connecting the three control volume centres considered for the application of the method, constitutes  a stylized "L".
Here, we follow  the description of the L-method given  in \cite{aavatsmark2002introduction}.\\
\\

\begin{itemize}
	
\item[] Let us consider the triangle $x_1x_2x_3$ (see Figure \ref{fig:triangle}), and a linear function $g$ defined over this triangle. we define

\begin{equation}
\label{xgrad}
\mathbf{X}\nabla g =\begin{bmatrix}
g(x_2)-g(x_1) \\
\\
g(x_3)-g(x_1)
\end{bmatrix}
\end{equation}

where

\begin{equation}
\mathbf{X}=\begin{bmatrix}
\big(x_2-x_1\big)^T\\
\\
\big(x_3-x_1\big)^T
\end{bmatrix}
\end{equation}

Thereby, the gradient of the linear function $g$ may be expressed as follows:
\begin{equation}
\label{grad}
\nabla g = \frac{1}{T} \Bigg(\nu_2\Big(g(x_2)-g(x_1)\Big)+\nu_3\Big(g(x_3)-g(x_1)\Big)\Bigg)
\end{equation} 

where $\nu_2,\nu_3$ are respectively the normal vector to $x_2-x_1$ and $x_3-x_1$ defined by

\begin{equation*}
\nu_2=\mathbf{R}\big(x_2-x_1\big),~~~~~~~~~\nu_3=-\mathbf{R}\big(x_3-x_1\big)
\end{equation*}

with

\begin{equation*}
\mathbf{R}=\begin{bmatrix}
0   &  &   &  &  1 \\
-1 &   &  &  & 0
\end{bmatrix}
\end{equation*}

and 

\begin{equation*}
T=\nu_2^T\mathbf{R}\nu_3
\end{equation*}
Let's notice that the matrix $\mathbf{R}$ is a rotation of angle $-\frac{\pi}{2}$. Thereby the vector $\nu_2$ and $\nu_3$ have same length with the vectors $x_2-x_1$ and $x_3-x_1$.

\begin{figure}[!h]
\begin{center}	
\begin{tikzpicture}
	\draw[blue,thick] (1,1)--(3,5)--(6,2)--(1,1);
	\draw[blue,thick,->] (2,3)--(4,2);
	\draw[blue,thick,->] (3.5,1.5)--(3,4);
	\node[above] at (0.8,0.6){${x}_1$};
	\node[above] at (6.2,2){${x}_2$};
	\node[above] at (3,5){$x_3$};
	\node[above] at (3.3,3.7){$\nu_2$};
	\node[above] at (4,2){$\nu_3$};
\end{tikzpicture}	
\end{center}
\caption{Triangle $x_1x_2x_3$}
\label{fig:triangle}
\end{figure}
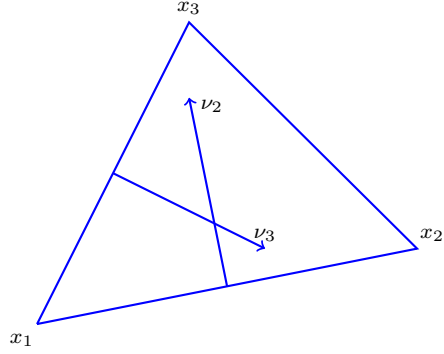
\end{itemize}
Let us called \textbf{interaction volume} $\R_{ij}$ a cell grid defined as follows
\begin{equation}
for~~~i,j=1,\ldots,N+1~~~~~\R_{ij}=[x_{i-1};x_i]\times[y_{j-1};y_j]
\end{equation}
We denote respectively by $x_1(x_{i-1},y_{j-1}),x_2(x_{i},y_{j-1}),x_3(x_{i},y_{j})$
~and~$ x_4(x_{i-1},y_{j})$ the centre of the control volume $\C_{ij},\C_{i+1,j},\C_{i,j+1}$ and $\C_{i+1,j+1}$.We denote also by $\bar{x}_1,\bar{x}_2,\bar{x}_3~and~\bar{x}_4$ the midpoints of the segment $x_1x_2$, $x_3x_4,x_1x_3~and~x_2x_4$.
We may notice that an  interaction volume $\R_{ij}$, for $i,j=1,\ldots,N+1$, is covering an area  in the intersection of the control volumes $\C_{ij},\C_{i+1,j},\C_{i,j+1}$ and $\C_{i+1,j+1}$. An interaction volume can be divided into 2 triangles such that the half edges $1,2$ are in the triangle $T_1=x_1x_2x_3$  and the half edges $3,4$ are in the triangle $T_2=x_1x_3x_4$ (see Figure \ref{fig:intervol}).\\
Here, We follow the procedure in \cite{aavatsmark2007multipoint}.
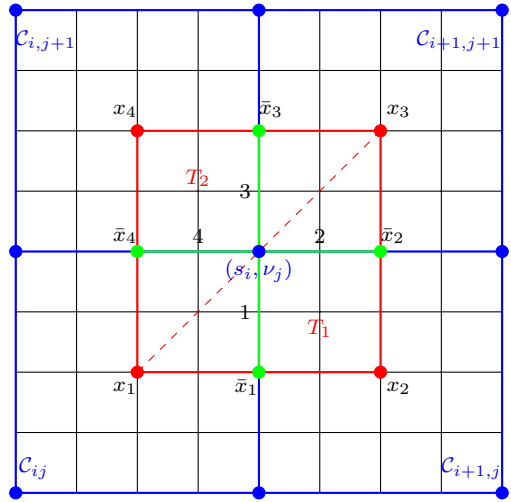
\begin{figure}[!h]
	\begin{center}
		\begin{tikzpicture}[scale=0.8]
		\draw (0,0) grid (8,8);
		\draw[blue,thick] (4,0)--(4,8);
		\draw[blue,thick] (0,4)--(8,4);
		\node[above,blue] at (0.3,0.1){$\mathcal{C}_{ij}$};
		\node[above,blue] at (7.5,0.1){$\mathcal{C}_{i+1,j}$};
		\node[below,blue] at (7.3,7.8){$\mathcal{C}_{i+1,j+1}$};
		\node[below,blue] at (0.5,7.8){$\mathcal{C}_{i,j+1}$};
		\node[below,red] at (5,3){$T_1$};
		\node[below,red] at (3,5.5){$T_2$};
		\draw[red,thick] (2,2)--(6,2);
		\draw[red,thick] (2,6)--(6,6);
		\draw[red,thick] (2,2)--(2,6);
		\draw[red,thick] (6,2)--(6,6);
		\draw[red,dashed] (2,2)--(6,6);
		\draw[blue,thick] (0,0)--(0,8);
		\draw[blue,thick] (0,0)--(8,0);
		\draw[blue,thick] (8,0)--(8,8);
		\draw[blue,thick] (0,8)--(8,8);
		\draw[green,thick] (2,4)--(6,4);
		\draw[green,thick] (4,2)--(4,6) ;
		\draw[red,fill=red] (6,6) circle (0.1);
		\draw[red,fill=red] (2,6) circle (0.1);
		\draw[red,fill=red] (6,2) circle (0.1);
		\draw[red,fill=red] (2,2) circle (0.1);
		\draw[blue,fill=blue] (4,4) circle (0.1)node[below]{$(s_{i},\nu_{j})$};
		\draw[blue,fill=blue] (0,0) circle (0.1);
		\draw[blue,fill=blue] (0,4) circle (0.1);
		\draw[blue,fill=blue] (0,8) circle (0.1);
		\draw[blue,fill=blue] (4,0) circle (0.1);
		\draw[blue,fill=blue] (4,8) circle (0.1);
		\draw[blue,fill=blue] (8,0) circle (0.1);
		\draw[blue,fill=blue] (8,4) circle (0.1);
		\draw[blue,fill=blue] (8,8) circle (0.1);
		\node[above] at (3,4){$4$};
		\node[left] at (4,5){$3$};
		\node[left] at (4,3){$1$};
		\node[above] at (5,4){$2$};
		\draw[green,fill=green] (4,6) circle (0.1);
		\draw[green,fill=green] (6,4) circle (0.1);
		\draw[green,fill=green] (4,2) circle (0.1);
		\draw[green,fill=green] (2,4) circle (0.1);
		\node[above] at (4.2,6.1){$\bar{x}_{3}$};
		\node[above] at (1.8,4){$\bar{x}_{4}$};
		\node[below] at (3.8,2){$\bar{x}_{1}$};
		\node[above] at (6.2,4){$\bar{x}_{2}$};
		\node[above] at (1.8,1.5){$x_1$};
		\node[above] at (6.3,1.5){$x_2$};
		\node[above] at (1.8,6.1){$x_4$};
		\node[above] at (6.3,6.1){$x_3$};
		\end{tikzpicture}
	\end{center}
	\caption{Interaction volume}
	\label{fig:intervol}
\end{figure}

In an interaction volume, we aim to compute the flux through the half edges $1,2,3$ and $4$ (See \ref{fig:intervol}).

Thereby, using (\ref{diffterm-disc}), the flux $f_p^{ij}$ through the half edge $p$ seen from the centre of the control volume $\C_{ij}$ is expressed as follows:

\begin{equation}
\label{flux}
f_p^{ij}=n_p^T\mathbf{M}^{ij}\nabla \V_{ij}
\end{equation}

where $n_p$ is the vector normal to the half edge $p$ with the same length.\\

Let us consider the triangle $T_1=x_1x_2x_3$ from the interaction volume $\R_{ij}$

\begin{itemize}

\item[]

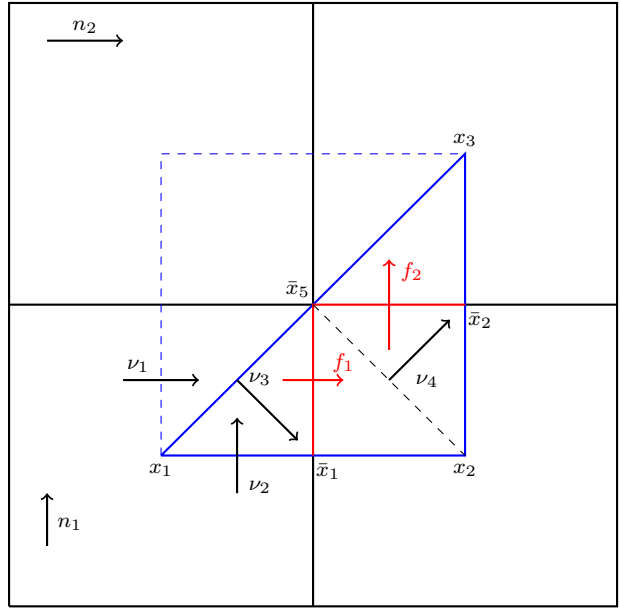
\begin{figure}[!h]
	\centering
	\begin{tikzpicture}
	\draw[black,thick] (0,0)--(8,0)--(8,8)--(0,8)--(0,0);
	\draw[black,thick] (0,4)--(4,4);
	\draw[black,thick] (6,4)--(8,4);
	\draw[black,thick] (4,0)--(4,2);
	\draw[black,thick] (4,4)--(4,8);
	\draw[blue,thick] (2,2)--(6,2)--(6,6)--(2,2);
	\draw[blue,dashed] (2,2)--(2,6)--(6,6);
	\draw[red,thick] (4,2)--(4,4)--(6,4);
	\draw[->,thick] (1.5,3)--(2.5,3);
	\draw[->,thick] (3,1.5)--(3,2.5);
	\draw[->,thick] (3,3)--(3.8,2.2);
	\draw[black,dashed] (4,4)--(6,2);
	\draw[->,thick] (5,3)--(5.8,3.8);
	\node[below] at (2,2){$x_1$};
	\node[below] at (6,2){$x_2$};
	\node[above] at (6,6){$x_3$};
	\node[above] at (3.8,4){$\bar{x}_5$};
	\node[below] at (4.2,2){$\bar{x}_1$};
	\node[below] at (6.2,4){$\bar{x}_2$};
	\node[above] at (1.7,3){$\nu_1$};
	\node[above] at (3.3,1.4){$\nu_2$};
	\node[above] at (5.5,2.8){$\nu_4$};
	\node[below] at (3.3,3.2){$\nu_3$};
	\draw[red,->,thick] (3.6,3)--(4.4,3);
	\draw[red,->,thick] (5,3.4)--(5,4.6);
	\node[above,red] at (4.4,3){$f_1$};
	\node[above,red] at (5.3,4.2){$f_2$};
	\draw[->,thick] (0.5,0.8)--(0.5,1.5);
	\draw[->,thick] (0.5,7.5)--(1.5,7.5);
	\node[above] at (0.8,0.9){$n_1$};
	\node[above] at (1,7.5){$n_2$};
	\end{tikzpicture}
	\caption{Triangle $T_1$}
	\label{fig:triT1}
\end{figure}

\item[] In the triangle $T_{12}=\bar{x}_1x_2\bar{x}_2$, using the flux expression (\ref{flux}) and following the gradient expression in (\ref{grad}) we get
\begin{eqnarray}
\label{fluxx212}
f_1^{i,j-1}&  =  & \omega_{12}^{i,j-1}\Big(\bar{\V}_2-\V_{i,j-1}\Big)-\omega_{11}^{i,j-1}\Big(\bar{\V}_1-\V_{i,j-1}\Big) \nonumber   \\
&    & \\
f_2^{i,j-1} & = &  \omega_{22}^{i,j-1}\Big(\bar{\V}_2-\V_{i,j-1}\Big)-\omega_{21}^{i,j-1}\Big(\bar{\V}_1-\V_{i,j-1}\Big) \nonumber
\end{eqnarray}
where

\begin{eqnarray*}
\omega_{11}^{i,j-1}= \frac{1}{T_1^{i,j-1}}\times n_1^T \mathbf{M}^{i,j-1} \nu_1 &~~~~~~~~~~~~~~~~~~~~~ &\omega_{12}^{i,j-1}= \frac{1}{T_1^{i,j-1}}\times n_1^T \mathbf{M}^{i,j-1} \nu_2  \\
&     & \\
&     & \\
\omega_{21}^{i,j-1}= \frac{1}{T_1^{i,j-1}}\times n_2^T \mathbf{M}^{i,j-1} \nu_1&~~~~~~~~~~~~~~~~~~~~~&\omega_{22}^{i,j-1}= \frac{1}{T_1^{i,j-1}}\times n_2^T \mathbf{M}^{i,j-1} \nu_2
\end{eqnarray*}

with

\begin{equation*}
T_1^{i,j-1}=-\nu_2^T\mathbf{R}\nu_1
\end{equation*}
 \item [] Moreover, using the property of the matrix $\mathbf{R}$, we have 
 
 \begin{equation}
 \nu_4=\mathbf{R}(\bar{x}_5-x_2)
 \end{equation}

Using also the equation (\ref{xgrad}) and the expression of gradient (\ref{grad}) in the triangle $T_{12}$ lead to 

\begin{equation}
\label{u5}
\bar{\V}_5 = \V_{i,j-1}+\chi_{42}^{i,j-1}\Big(\bar{\V}_2-\V_{i,j-1}\Big)-\chi_{41}^{i,j-1}\Big(\bar{\V}_1-\V_{i,j-1}\Big)
\end{equation}

where

\begin{eqnarray}
\chi_{41}^{i,j-1}=\frac{1}{T_1^{i,j-1}}\nu_4^T\mathbf{R}\nu_1~~~~~~~~~~~~~~~~~~~~~~\chi_{42}^{i,j-1}=\frac{1}{T_1^{i,j-1}}\nu_4^T\mathbf{R}\nu_2 \nonumber
\end{eqnarray}

 \item [] In the triangle $T_{11}=x_1\bar{x}_1\bar{x}_5$, we have

\begin{eqnarray}
\label{fluxx1}
f_1^{i-1,j-1}   =   \omega_{13}^{i-1,j-1}\Big(\bar{\V}_1-\V_{i-1,j-1}\Big)+\omega_{12}^{i-1,j-1}\Big(\bar{\V}_5-\V_{i-1,j-1}\Big) 
\end{eqnarray}

Replacing $\bar{\V}_5$ by its expression (\ref{u5}) in  (\ref{fluxx1}),we have

\begin{eqnarray}
\label{fluxx11}
f_1^{i-1,j-1}   =   \omega_{13}^{i-1,j-1}\Big(\bar{\V}_1-\V_{i-1,j-1}\Big)+\omega_{12}^{i-1,j-1}\Bigg(\V_{i,j-1}-\V_{i-1,j-1}+\chi_{42}^{i,j-1}
\Big(\bar{\V}_2-\V_{i,j-1}\Big)-\chi_{41}^{i,j-1}\Big(\bar{\V}_1-\V_{i,j-1}\Big)\Bigg)\nonumber\\
\end{eqnarray}

where

\begin{eqnarray*}
	\omega_{13}^{i-1,j-1}= \frac{1}{T_1^{i-1,j-1}}\times n_1^T \mathbf{M}^{i-1,j-1} \nu_3 &~~~~~~~~~~~~~~~~~~~~~ &\omega_{12}^{i-1,j-1}= \frac{1}{T_1^{i-1,j-1}}\times n_1^T \mathbf{M}^{i-1,j-1} \nu_2  \\
\end{eqnarray*}

with

\begin{equation*}
T_1^{i-1,j-1}=\nu_3^T\mathbf{R}\nu_2
\end{equation*}

\item [] Similarly, in the triangle $T_{13}=\bar{x}_5\bar{x}_2x_3$

\begin{eqnarray}
\label{fluxx2}
f_2^{ij}   =   -\omega_{21}^{ij}\Big(\bar{\V}_5-\V_{ij}\Big)+\omega_{23}^{ij}\Big(\bar{\V}_2-\V_{ij}\Big) 
\end{eqnarray}

Replacing $\bar{\V}_5$ by its expression (\ref{u5}) in  (\ref{fluxx2}),we have

\begin{eqnarray}
\label{fluxx22}
f_2^{ij}   =   -\omega_{21}^{ij}\Bigg(\V_{i,j-1}-\V_{ij}+\chi_{42}^{i,j-1}\Big(\bar{\V}_2-\V_{i,j-1}\Big)-\chi_{41}^{i,j-1}
\Big(\bar{\V}_1-\V_{i,j-1}\Big)\Bigg)+\omega_{23}^{ij}\Big(\bar{\V}_2-\V_{ij}\Big)\nonumber\\
\end{eqnarray}

where

\begin{eqnarray*}
	\omega_{21}^{ij}= \frac{1}{T_1^{ij}}\times n_2^T \mathbf{M}^{ij} \nu_1 &~~~~~~~~~~~~~~~~~~~~~ &\omega_{23}^{ij}= \frac{1}{T_1^{ij}}\times n_2^T \mathbf{M}^{ij} \nu_3  \\
\end{eqnarray*}

with

\begin{equation*}
T_1^{ij}=-\nu_3^T\mathbf{R}\nu_3
\end{equation*}

\item [] Since the flux is continuous through edges, then using (\ref{fluxx212}), (\ref{fluxx11}) and (\ref{fluxx22}) we have

\begin{eqnarray}
\label{sys-eq1}
\left\lbrace
\begin{array}{l}
f_1   =  \omega_{12}^{i,j-1}\Big(\bar{\V}_2-\V_{i,j-1}\Big)-\omega_{11}^{i,j-1}\Big(\bar{\V}_1-\V_{i,j-1}\Big)
   \\
    ~~~~=   \omega_{13}^{i-1,j-1}\Big(\bar{\V}_1-\V_{i-1,j-1}\Big)+\omega_{12}^{i-1,j-1}\Bigg(\V_{i,j-1}-\V_{i-1,j-1}
    +\chi_{42}^{i,j-1}\Big(\bar{\V}_2-\V_{i,j-1}\Big)-\chi_{41}^{i,j-1}\Big(\bar{\V}_1-\V_{i,j-1}\Big)\Bigg)\ \\
            \\
             \\
f_2  =  \omega_{22}^{i,j-1}\Big(\bar{\V}_2-\V_{i,j-1}\Big)-\omega_{21}^{i,j-1}\Big(\bar{\V}_1-\V_{i,j-1}\Big) \\
     \\
 ~~~~=   -\omega_{21}^{ij}\Bigg(\V_{i,j-1}-\V_{ij}+\chi_{42}^{i,j-1}\Big(\bar{\V}_2-\V_{i,j-1}\Big)-\chi_{41}^{i,j-1}\Big(\bar{\V}_1-\V_{i,j-1}
 \Big)\Bigg)+\omega_{23}^{ij}\Big(\bar{\V}_2-\V_{ij}\Big)
\end{array}\right.
\end{eqnarray}

\end{itemize}

by setting 

\begin{equation}
g=\left[\begin{array}{c} f_1 \\ f_2\end{array}\right],~~~~~~~~~~W=\left[\begin{array}{c} \V_{i-1,j-1} \\   \V_{i,j-1} \\  \V_{ij}\end{array}\right],~~~~~~~~~~~~~~~~ \V=\left[\begin{array}{c} \bar{\V}_1 \\  \bar{\V}_2 \\  \end{array}\right]
\end{equation}

The system of equations (\ref{sys-eq1}) can be written as

\begin{equation}
\label{syseq-flux}
g=C^{ij} \V+D^{ij}W
\end{equation}

where

\begin{eqnarray*}
		C^{ij} & = & \left[\begin{array}{cc}
			-\omega_{11}^{i,j-1} & \omega_{12}^{i,j-1} \\
			&    \\
			-\omega_{21}^{i,j-1}  & \omega_{22}^{i,j-1}\\
		\end{array}\right]~~~~~~~~~~~~~~~~~~	D^{ij}=\left[\begin{array}{ccc}
		0 & \omega_{11}^{i,j-1}-\omega_{12}^{i,j-1} & 0  \\
		&   &   \\
		0 &  \omega_{21}^{i,j-1}-\omega_{22}^{i,j-1} & 0 \\
		  \end{array}\right]
\end{eqnarray*}

Using the expressions at both sides of the second equalities of system equations (\ref{sys-eq1}) we have

\begin{equation}
\label{syseq-mat}
A^{ij}\V=B^{ij}W
\end{equation}

\begin{eqnarray*}
	A^{ij} & = & \left[\begin{array}{cc}
		\omega_{11}^{i,j-1}+\omega_{13}^{i-1,j-1}-\omega_{12}^{i-1,j-1}\chi_{41}^{i,j-1} & -\omega_{12}^{i,j-1}+\omega_{12}^{i-1,j-1}\chi_{42}^{i,j-1} 	 \\
		&    \\
	\omega_{21}^{i,j-1}+\omega_{21}^{ij}\chi_{41}^{i,j-1}	  & -\omega_{22}^{i,j-1}+\omega_{23}^{ij}-\omega_{21}^{ij}\chi_{42}^{i,j-1}\\
	\end{array}\right]
\end{eqnarray*}

\begin{eqnarray*}
	B^{ij} & = & \left[\begin{array}{cccc}
	\omega_{13}^{i-1,j-1}+\omega_{12}^{i-1,j-1} &   & -\omega_{12}^{i,j-1}+\omega_{11}^{i,j-1}-\omega_{12}^{i-1,j-1}\big(1+\chi_{41}^{i,j-1}-\chi_{42}^{i,j-1}\big) &  0 	 \\
		&    &    & \\
		0 &   & -\omega_{22}^{i,j-1}+\omega_{21}^{i,j-1}+\omega_{21}^{ij}\big(1+\chi_{41}^{i,j-1}-\chi_{42}^{i,j-1}\big)  	  &   -\omega_{21}^{ij}+\omega_{23}^{ij}\\
	\end{array}\right]
\end{eqnarray*}

Thereby, by solving (\ref{syseq-mat}) with respect to $\V$ and replacing in (\ref{syseq-flux}) we get

\begin{equation}
\label{eq-trans1}
g=R^{ij}\V
\end{equation}

where

\begin{equation*}
R^{ij}=C^{ij}[A^{ij}]^{-1}B^{ij}+D^{ij}
\end{equation*}

Now considering the triangle $T_2$ (see figure \ref{fig:intervol}) and applying the above procedure used in the triangle $T_1$, we are able to the fluxes through the half edges  $3$ and $4$ as follows:

\begin{equation}
\label{eq-trans2}
h=S^{ij}V
\end{equation}

where

\begin{equation*}
h=\left[\begin{array}{c} f_3 \\ f_4\end{array}\right],~~~~~~~~~~V=\left[\begin{array}{c} \V_{i,j-1} \\   \V_{ij} \\  \V_{ij}\end{array}\right]
\end{equation*}

For simplicity, in an interaction volume $\R_{ij}$, the flux through the half edges $1,2,3$ and $4$ are given by

\begin{equation}
\label{eq-trans}
f=T^{ij}\V
\end{equation}
where

\begin{equation}
f=\left[\begin{array}{c} f_1 \\ f_2 \\ f_3 \\ f_4 \end{array}\right],~~~~~~~~~~\V=\left[\begin{array}{c} \V_{i-1,j-1} \\   \V_{i,j-1} \\ \V_{ij} \\ \V_{i-1,j} \end{array}\right]
\end{equation}

 and $T^{ij}$ is $4\times 4$ matrix coming from $R^{ij}$ and $S^{ij}$ defined in (\ref{eq-trans1}),(\ref{eq-trans2}).
 $T^{ij}$ is called the transmissibility matrix of the interaction volume $\R_{ij}$.\\

Let us notice that the flux through a full  edge will be the addition of the fluxes through its 2 half edges.\\
\\
Let us recall that, from (\ref{diffterm-disc}), our aim is to compute the flux through the edges of the control volume $\C_{ij}$. In order to cover the boundary of a control volume, we need to four interaction volumes

\begin{figure}[!h]
	\centering
	\begin{tikzpicture}[scale=0.5]
	\draw[black,thick] (0,0)--(0,16)--(16,16)--(16,0)--(0,0);
	\draw[black,thick] (4,0)--(4,16);
	\draw[black,thick] (8,0)--(8,16);
	\draw[black,thick] (12,0)--(12,16);
	\draw[black,thick] (0,4)--(16,4);
	\draw[black,thick] (0,8)--(16,8);
	\draw[black,thick] (0,12)--(16,12);
	\draw[red,thick] (4,4)--(4,12)--(12,12)--(12,4)--(4,4);
	\draw[red,thick] (8,4)--(8,12);
	\draw[red,thick] (4,8)--(12,8);
	\draw[blue,dotted] (10,6)--(6,6)--(6,10)--(10,10);
	\draw[black,thick]  (10,6)--(10,10);
	\node[below,black] at (10.2,8){$\E$};
	\draw[red,fill=red] (4,4) circle (0.1);
	\draw[red,fill=red] (4,8) circle (0.1);
	\draw[red,fill=red] (4,12) circle (0.1);
	\draw[red,fill=red] (8,4) circle (0.1);
	\draw[red,fill=red] (8,8) circle (0.1);
	\draw[red,fill=red] (8,12) circle (0.1);
	\draw[red,fill=red] (12,4) circle (0.1);
	\draw[red,fill=red] (12,8) circle (0.1);
	\draw[red,fill=red] (12,12) circle (0.1);
	\node[below,red] at (8.7,8){$(x_i,y_j)$};
	\node[above,blue] at (6.7,6){$\C_{ij}$};
	\node[above] at (9.8,6.2){$3$};
	\node[above] at (9.8,8.2){$1$};
	\draw[black,fill=black] (10,8) circle (0.1);
	\node[above,red] at (4.5,4.1){$\R_{ij}$};
	\node[above,red] at (13.5,4.3){$\R_{i+1,j}$};
	\node[above,red] at (13.7,8.7){$\R_{i+1,j+1}$};
	\node[above,red] at (4.5,8.7){$\R_{i,j+1}$};
	\end{tikzpicture}
	\caption{}
\end{figure}

 Let us denote, for the volume control $\C_{ij}$, by ${}_{\E}f_{l}^{ij}$ the flux through lower half eastern edge, by ${}_{\E}f_{u}^{ij}$ the flux through the upper  half eastern edge. The flux  ${}_{\E}f^{ij}$  through the east edge of the control volume $\C_{ij}$ is calculated as follows:

 The lower half eastern edge is in position 3 in the triangle $T_2$  of the interaction volume $\R_{i+1,j}$  (See figure 4).
 So by using \eqref{eq-trans} we have:

 \begin{equation*}
 {}_{\E}f_{l}^{ij}=T_{31}^{i+1,j}\V_{i,j-1}+T_{33}^{i+1,j}\V_{i+1,j}+T_{34}^{i+1,j}\V_{ij}
 \end{equation*}
 
 Similarly, the upper half eastern edge  is  in position 1 in the triangle $T_1$  of the interaction volume $\R_{i+1,j+1}$  and it is in position 1 in the interaction volume. Thereby using (\eqref{eq-trans}) we have: 
 
 \begin{equation*}
 {}_{\E}f_u^{ij}=T_{11}^{i+1,j+1}\V_{ij}+T_{12}^{i+1,j+1}\V_{i+1,j}+T_{13}^{i+1,j+1}\V_{i+1,j+1}
 \end{equation*}
 Finally the flux through the eastern edge of the control volume $\C_{ij}$ will be the addition of ${}_{\E}f_{l}^{ij}$ and  ${}_{\E}f_u^{ij}$. Thereby we have:

\begin{eqnarray}
{}_{\E}f^{ij} & = &  {}_{\E}f_d^{ij}+{}_{\E}f_u^{ij}\nonumber \\
&   & \nonumber \\
{}_{\E}f^{ij} & = &  T_{31}^{i+1,j}\V_{i,j-1}+\Big(T_{33}^{i+1,j}+T^{i+1,j+1}_{12}\Big)\V_{i+1,j}+\Big(T_{11}^{i+1,j+1}+T_{34}^{i+1,j}\Big)\V_{ij}
+T_{13}^{i+1,j+1}\V_{i+1,j+1} \nonumber \\
\end{eqnarray}

The same method is applied to calculate the flux through the northern, western and southern  edges  of the control volume $\C_{ij}$. The flux through the edges of the control volume $\C_{ij}$ is obtained by summing up  the flux through the 4 edges. This gives :

\begin{equation}
\label{mpfa-L}
\F_{ij}  = a_{ij}\V_{ij}+b_{ij}\V_{i+1,j}+c_{ij}\V_{i+1,j+1}+d_{ij}\V_{i,j+1}+e_{ij}\V_{i-1,j}+\alpha_{ij}\V_{i-1,j-1}+\beta_{ij}\V_{i,j-1} 
\end{equation}

where

\begin{eqnarray*}
	&    & a_{ij}=T_{11}^{i+1,j+1}+T_{34}^{i+1,j}+T_{41}^{i+1,j+1}+T_{22}^{i,j+1}-T_{12}^{i,j+1}-T_{33}^{ij}-T_{23}^{ij}-T_{44}^{i+1,j}\\
	&    & \\
	&    & b_{ij}=T_{33}^{i+1,j}+T^{i+1,j+1}_{12}-T_{43}^{i+1,j}~~~~~~~~~c_{ij}=T_{13}^{i+1,j+1}+T^{i+1,j+1}_{43}~~~~~~~~d_{ij}=T_{23}^{i,j+1}+T_{44}^{i+1,j+1}-T_{13}^{i,j+1}\\
	&     & \\
	&     & e_{ij}=T_{21}^{i,j+1}-T_{11}^{i,j+1}-T_{34}^{ij}~~~~~~~~\alpha_{ij}=-T_{31}^{ij}-T_{21}^{ij}~~~~~~~~\beta_{ij}=T_{31}^{i+1,j}-T_{22}^{ij}-T_{41}^{i+1,j}
\end{eqnarray*}

This leads to  a system of equations which can be written as follows:

\begin{equation}
\label{eq-MPFA-L}
\F=A_{mp}\V+F_{mp}
\end{equation}

where 
\begin{equation*}
\F=\begin{bmatrix}
\F_{11}\\
\F_{12}\\
\vdots\\
\F_{1N}\\
\F_{21}\\
\F_{22}\\
\F_{NN}
\end{bmatrix}~~
\V=\begin{bmatrix}
\V_{11}\\
\V_{12}\\
\vdots\\
\V_{1N}\\
\V_{21}\\
\V_{22}\\
\V_{NN}
\end{bmatrix}~~
A_{mp}=\begin{bmatrix}
W_1 & X_1 & 0_N & \ldots & \ldots & \ldots & \ldots & 0_N\\ 
Y_2 & W_2 & X_2 & \ddots &  &  &  & \vdots  \\ 
0_N & Y_3 & W_3 & X_3 & \ddots &  & & \vdots\\ 
\vdots & \ddots & Y_4 & W_4 & X_4 & \ddots & & \vdots\\ 
\vdots &  & \ddots & \ddots & \ddots & \ddots & \ddots & \vdots \\ 
\vdots &  &  & \ddots& \ddots & \ddots & \ddots & 0_N&  \\ 
\vdots &  &  &  & \ddots & Y_{N-1} & W_{N-1} & X_{N-1} \\ 
0_N & \ldots & \ldots & \ldots & \ldots & 0_N & Y_N & W_N
\end{bmatrix}  
\end{equation*}

with $0_N$ is $N\times N $ zeros matrix, $W_i,X_i,Z_i$ are tridiagonal matrix  and $F_{mp}$ is a $N^2$ vector coming from the boundary conditions.\\
The diffusion matrix $A_{mp}$ is under the folowing form:

\begin{figure}[hbtp]
	 \centering
\includegraphics[scale=0.5]{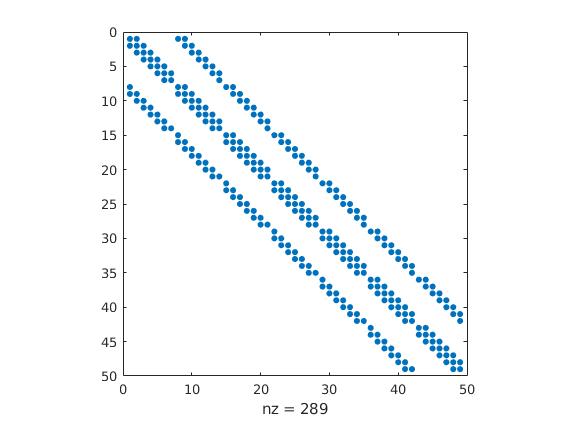}
  \caption{A structure of the  diffusion matrix using L-MPFA method.}
  \label{spy-LMPFA}
\end{figure}

As we can see on Fig.\ref{spy-LMPFA}, the L-MPFA method is a 7 points stencil method, unlikely to the O-MPFA method (see \cite{rock2019mpfa} ) which is a 9 points stencil method.


\subsection{Discretisation of the convection term}

The integral of  convection term  

\begin{equation*}
\int_{\C_{ij}}\nabla (f \V)d\C
\end{equation*}

where 
\begin{equation*}
f=\left(\begin{array}{c}
(r-\sigma_1^2-\frac{1}{2}\rho\sigma_1\sigma_2)x \\ \\ (r-\sigma_2^2-\frac{1}{2}\rho\sigma_1\sigma_2)y
\end{array}\right)
\end{equation*}

will be approximated  using the upwind methods ($1^{st}$ and $2^{nd}$ order ). We start by applying the  divergence theorem, and we have  for $i,j=1,...,N$:

\begin{equation}
\label{div-theo-up}
I^{ij}=\int_{\C_{ij}}\nabla (f \V) d \C_{ij}=\int_{\partial \C_{ij}}(f\cdot\V)\cdot\vec{n}d\partial \C_{i,j}
\end{equation}

with  \vec{n} an outward unit normal vector.

\subsubsection{First order upwind method}
 
The \textbf{first order upwind method} discussed in \cite[chapter 4.8]{leveque2004finite} will be applied to evaluate the second term of \eqref{eqfinvol}.\\
\\
$I^{ij}$ is calculated by summing up the flux through the edges of the control volume $\C_{ij}$.\\
The flux through an edge using the first order upwind will depend on the sign of $f\cdot\vec{n}$ on this edge. If the sign of $f\cdot\vec{n}$ is positive, $\V_{ij}$ will be used to approximate $(f\cdot \vec{n} \V)$ otherwise we will use the value of $\V$ in other side of the edge.\\

By doing so, we get $\forall i,j=1,...,N$
\begin{eqnarray}
\label{flux-up1}
I^{ij} & =  & \epsilon_{ij}\V_{i-1,j}+\mu_{ij}\V_{i,j-1}+\Omega_{ij}\V_{ij}+
\phi_{ij}\V_{i,j+1}+\Psi_{ij}\V_{i+1,j}
\end{eqnarray}

where 
\begin{eqnarray*}
	&   &\epsilon_{ij}=-l_jf_x^{i-1}\max(f_x^{i-1},0);~~~~~\mu_{ij}=-h_if_y^{j-1}\max(f_y^{j-1},0)\\
	&   & \\
	&   & \Omega_{ij}=l_j\Bigg(f_x^{i}\max(f_x^{i},0)- f_x^{i-1}\min(f_x^{i-1},0)\Bigg)
	+h_i\Bigg(f_y^j\max(f_y^j,0)-f_y^{j-1}\min(f_y^{j-1},0)\Bigg)\\
	&  & \\
	&   & \\
	&   & \phi_{ij}=h_if_y^j\min(f_y^j,0);~~~~~~~~
	\Psi_{ij}=l_jf_x^{i}\min(f_x^{i},0)
\end{eqnarray*}


Equation \eqref{flux-up1} will lead to a system of equations which be written as follows:

\begin{equation}
\label{flux-up1-mat}
I=A_{up} \V+F_{up}
\end{equation}

where  $A_{up}$ is a $M\times N$ matrix 

\begin{equation*}
I=\begin{bmatrix}
I^{11}\\
I^{12}\\
\vdots\\
I^{1N}\\
I^{21}\\
I^{22}\\
\vdots\\
\vdots\\
I^{NN}\\
\end{bmatrix}
~~~~~~\V=\begin{bmatrix}
\V_{11}\\
\V_{12}\\
\vdots\\
\V_{1N}\\
\V_{21}\\
\V_{22}\\
\vdots\\
\vdots\\
\V_{NN}\\
\end{bmatrix}A_{up}=\begin{bmatrix}
H_1 & P_1 & 0_N & \ldots & \ldots & \ldots & 0_N \\ 
Q_2 & H_2 & P_2& \ddots & &   & \vdots \\
0_N & Q_3 & H_3 & P_3 & \ddots &   &   \vdots \\
\vdots & \ddots & \ddots & \ddots & \ddots & \ddots & \vdots \\
\vdots  &   & \ddots & Q_{N-2}  & H_{N-2} & P_{N-2} & 0_N \\
\vdots &  &  & \ddots & Q_{N-1} & H_{N-1} & P_{N-1}\\
0_N & \dots   &   \ldots  & \ldots  &  0_N &  Q_N  & H_N
\end{bmatrix}
\end{equation*}
\\
with $0_{N}$ is $N\times N$ zeros matrix, $H_i$ is a tridiagonal matrix, $P_i,Q_i$ are diagonal matrices and $F_{up}$ is a vector coming from the boundary conditions. \\
The advection matrix $A_{up}$ using the first order upwind method is under following form

Therefore, combining the L-MPFA method \eqref{eq-MPFA-L} and the first order upwind \eqref{flux-up1-mat}, we get 

\begin{equation}
\label{mpfa-up1}
\frac{d\V}{d\tau}=A\V+G(\V)+F
\end{equation} 

with 

\begin{eqnarray*}
	A=L^{-1}\Bigg(A_{mp}+A_{up}+A_L\Bigg)~~~~G(\V)=\beta \Big[\max\Big(\V^*-\V,0\Big)\Big]^{1/k}~~F=L^{-1}\Bigg(F_{mp}+F_{up}\Bigg)
\end{eqnarray*}
\\
\\
where $A_L$ is a diagonal matrix of size $N\times N$ coming from the discretisation of \eqref{linearterm}. The diagonal elements of $A_L$ are $A_{ii}=h_il_i\lambda$ for $i=1,...,N$ with $\lambda$ given in \eqref{conservation}. The matrix L is also a diagonal matrix of size $N\times N$ whose diagonal elements are $L_{ii}=h_il_i$ for $i=1,\ldots,N$

\subsubsection{Second order upwind method}

A second order approximation is used to calculated  the flux defined in (\ref{div-theo-up}). For instance, the flux ${}_{\E}J^{ij}$ through the eastern edge of the control volume $\C_{ij}$ is approximated as follows:
We will start by  giving  an approximation of the gradient in the integral expression \eqref{diffterm-disc}.
\begin{equation}
{}_{\E}J^{ij}=\int_{\E_{ij}} (f\cdot \V)\cdot n_{\E}
\end{equation}

where

\begin{equation*}
\vec{n}_{\E}=\begin{bmatrix}
1\\
\\
0
\end{bmatrix}
\end{equation*}

is the outward unit normal vector to the eastern edge $\E_{ij}$ of the control volume $\C_{ij}$. We set $f_x=f\cdot n_{\E}$ and we have

\begin{equation}
\V \approx \left\lbrace 
\begin{array}{c}
\frac{3\V_{ij}-\V_{i-1,j}}{2}~~~~~if~~f_x\geq 0 \\
\\
\\
\frac{3\V_{i+1,j}-\V_{i+2,j}}{2}~~~~~~if~~f_x<0
\end{array} \right.
\end{equation}

Then we get

\begin{equation}
{}_{\E}J^{ij}=l_j\Bigg[\frac{3}{2}\max(f_x^{i+1},0)\V_{ij}-\frac{1}{2}\max(f_x^{i+1},0)\V_{i-1,j}
+\frac{3}{2}\min(f_x^{i+1},0)\V_{i+1,j}-\frac{1}{2}\min(f_x^{i+1},0)\V_{i+2,j}\Bigg]
\end{equation}

with
\begin{equation*}
f_x^{i+1}=(r-\sigma_x^2 - \frac{1}{2}\rho\sigma_x \sigma_y)x_{i+\frac{1}{2}}
\end{equation*}

We use the same argument to calculate the flux ${}_{\N}J^{ij},{}_{\W}J^{ij},{}_{S}J^{ij}$ through the  northern, western and southern edges of the control volume $\C_{ij}$ and after sum them up. We get then

\begin{eqnarray}
\label{eqflux-2ndup}
	J^{ij} &  =  & \epsilon_{ij}\V_{i-2,j}+\eta_{ij}\V{i-1,j}+\kappa_{ij}\V_{i,j-2}+\mu_{ij}\V_{i,j-1}+\Omega_{ij}\V_{ij}
+\phi_{ij}\V_{i,j+1}+\Psi_{ij}\V_{i,j+2}+\Delta_{ij}\V_{i+1,j} \nonumber\\
&    &  \nonumber \\
&    & +\Pi_{ij}\V_{i+2,j}
\end{eqnarray}

where

\begin{eqnarray*}
	&    & \epsilon_{ij}=\frac{1}{2}l_j\max(f_x^i,0)~~~~~~~~~~~~~\eta_{ij}=-\frac{1}{2}l_j\max(f_x^{i+1},0)-\frac{3}{2}l_j\max(f_x^{i},0)\\
	&    & \\
	&    & \\
	&    & \kappa_{ij}=\frac{1}{2}h_i\max(f_y^j,0)~~~~~~~~~~~~~~~~~\mu_{ij}=-\frac{1}{2}h_i\max(f_y^{j+1},0)-\frac{3}{2}h_i\max(f_y^{j},0)\\
	&    & \\
	&    & \\
	&    & \Omega_{ij}=\frac{3}{2}l_j\max(f_x^{i+1},0)+\frac{3}{2}h_i\max(f_y^{j+1},0)-\frac{3}{2}l_j\min(f_x^{i},0)-\frac{3}{2}h_i\min(f_y^{j},0)\\
	&    & \\
	&    & \\
	&     & \phi_{ij}=\frac{3}{2}h_i\min(f_y^{j+1},0)+\frac{1}{2}h_i\min(f_y^j,0)~~~~~~~~~~~~~~~~\Psi_{ij}=-\frac{1}{2}h_i\min(f_y^{j+1},0)\\
	&     & \\
	&     & \\
	&     & \Delta_{ij}=\frac{3}{2}l_j\min(f_x^{i+1},0)+\frac{1}{2}l_j\min(f_x^i,0)~~~~~~~~~~\Pi_{ij}=-\frac{1}{2}l_j\min(f_x^{i+1},0)
\end{eqnarray*}

For the control volumes near the boundary of the study domain, the first order upwind method is used for the approximation of the flux through edges directly connected to the boundary.

Equation (\ref{eqflux-2ndup}) leads to a system of equations which can be written as:

\begin{equation}
\label{flux-up2-mat}
J=A_{2up} \V+F_{2up}
\end{equation}

where 

\begin{equation*}
J=\begin{bmatrix}
J^{11}\\
J^{12}\\
\vdots\\
J^{1N}\\
J^{21}\\
J^{22}\\
\vdots\\
J^{NN}
\end{bmatrix}
~~~~~~
\V=\begin{bmatrix}
\V_{11}\\
\V_{12}\\
\vdots\\
\V_{1N}\\
\V_{21}\\
\V_{22}\\
\vdots\\
\V_{NN}\\
\end{bmatrix}~~~~
A_{2up}=\begin{bmatrix}
K_1 & R_1 & G_1 & 0_N & \ldots & \ldots & 0_N \\ 
S_2 & K_2 & R_2& G_2 & & \ddots  & \vdots  \\
H_3 & S_3 & K_3 & R_3 & G_3 &   &   \vdots\\
\vdots & \ddots & \ddots & \ddots & \ddots & \ddots & \vdots \\
\vdots  &   & H_{N-2} & S_{N-2}  & K_{N-2} & R_{N-2} & G_{N-2} \\
\vdots &  &  & H_{N-1} & S_{N-1} & K_{N-1} & R_{N-1}& \\
0_N & \dots   &   \ldots  & 0_N  &  H_N &  S_N  & K_N \\
\end{bmatrix}
\end{equation*}

where for $i=1,\ldots,N~K_i$ is penta-diagonal matrice and $R_i,G_i,S_i,H_i$ are diagonal matrices. $F_{2up}$ is a vector coming from the boundary conditions.

Therefore, combining the L-MPFA method  \eqref{eq-MPFA-L} and the second order upwind method \eqref{flux-up2-mat}, we have 

\begin{equation}
\label{mpfa-up2}
\frac{d\V}{d\tau}=A\V+G(\V)+F
\end{equation} 

with 

\begin{eqnarray*}
	A=L^{-1}\Bigg(A_{mp}+A_{2up}+A_L\Bigg)~~~~G(\V)=\beta \Big[\max\Big(\V^*-\V,0\Big)\Big]^{1/k}~~F=L^{-1}\Bigg(F_{mp}+F_{2up}\Bigg)
\end{eqnarray*}
\\
where $A_L$ is a diagonal matrix of size $N\times N$ coming from the discretisation of \eqref{linearterm}. The diagonal elements of $A_L$ are $A_{ii}=h_il_i\lambda$ for $i=1,...,N$ with $\lambda$ given in \eqref{conservation}. The matrix L is also a diagonal matrix of size $N\times N$ whose diagonal elements are $L_{ii}=h_il_i$ for $i=1,\ldots,N$
\\
\\
Besides, the ellipticity condition for the PDE (\ref{BS-op}) is not satisfied when the stocks price  ($x\rightarrow 0$ and/or $y\rightarrow 0$) is near to zero.  This may cause some oscillations of the numerical solution when the PDE is degenerate.\\
Nevertheless, \cite{wang2004novel} suggested a fitted finite volume method to deal with the degeneracy of the PDE. Thereby, the fitted finite volume method will be applied in the degeneracy region ($x\rightarrow 0$ and/or $y\rightarrow 0$) in the next section.

\subsubsection{Fitted finite volume }

The fitted finite volume method is used to approximated the flux through edges which are (fully) in the degeneracy region i.e the western edge of the control volume $\C_{1,j}~~~j=1,..,N$ and the southern edge of the control volume $\C_{i,1}~~~i=1,\ldots,N$.\\
\\
For the western edge of the control volume $\C_{1,j}~~~~j=1,\ldots,N$, using the mid-quadrature rule, we have:

\begin{eqnarray}
\label{int-wes}
	\int_{(x_{\frac{1}{2}},y_{j-\frac{1}{2}})}^{(x_{\frac{1}{2}},y_{j+\frac{1}{2}})}\Bigg(m_{11}\frac{\partial \V}{\partial x}+m_{12}\frac{\partial \V}{\partial y}+p\V\Bigg)dy \approx \Bigg(m_{11}\frac{\partial \V}{\partial x}+m_{12}\frac{\partial \V}{\partial y}+p\V\Bigg)_{(x_{\frac{1}{2}},y_j)} \cdot l_j
\end{eqnarray}

besides, we have:

\begin{eqnarray}
\label{intg-wes}
	m_{11}\frac{\partial \V}{\partial x}+m_{12}\frac{\partial \V}{\partial y}+p\V& = & x\Bigg(ax\frac{\partial \V}{\partial x}+d\frac{\partial \V}{\partial y}+b\V\Bigg)
\end{eqnarray}

	with $a=\frac{1}{2}\sigma_1^2$ , $b=r-\sigma_1^2-\frac{1}{2}\rho\sigma_1\sigma_2$ and $d=\frac{1}{2}\rho\sigma_1\sigma_2y$.\\
\\
We want to approximate

\begin{equation*}
g(\V)=ax\frac{\partial \V}{\partial x}+b\V
\end{equation*}

by a constant over $I_{x_1}=(0,x_1)$ satisfying the following two-point boundary value problem

\begin{align}
\label{bvp}
\left\lbrace \begin{array}{l}
g'(v)~~~~=\Bigg(ax\frac{\partial v}{\partial x}+bv\Bigg)'=K_1\\
\\
\\
v(0,y_j) =\V_{0,j}~~~~~~~~~v(x_1,y_j)=\V_{1,j}
\end{array}\right.
\end{align}

By solving this problem we get 

\begin{eqnarray}
\label{sol-bvp}	
	\V=\V_{0,j}+(\V_{1,j}-\V_{0,j})\frac{x}{x_{1}}
\end{eqnarray}

Thereby, by using (\ref{int-wes}),(\ref{intg-wes}), (\ref{bvp}), (\ref{sol-bvp}) and the forward difference to approximate the first partial derivative  $\frac{\partial \V}{\partial y}$ we get 

%

\begin{equation}
\label{west-appr-fitfin}
\int_{(x_{\frac{1}{2}},y_{j-\frac{1}{2}})}^{(x_{\frac{1}{2}},y_{j+\frac{1}{2}})}\Bigg(m_{11}\frac{\partial \V}{\partial x}+m_{12}\frac{\partial \V}{\partial y}+p\V\Bigg)dy \approx \frac{1}{2}x_1\Big[\frac{1}{2}l_j(a+b)-d_j\Big]\V_{i,1}+\frac{1}{2}d_jx_1\V_{1,j+1}-\frac{1}{4}x_1l_j(a-b)\V_{0,j}
\end{equation}

where 

\begin{eqnarray*}
a=\frac{1}{2}\sigma_1^2,~~~~~~~b=r-\sigma_1^2-\frac{1}{2}\rho\sigma_1\sigma_2~~~~~~~~d_j=\frac{1}{2}\rho\sigma_1\sigma_2y_j~~~~~~~l_j=y_{j+\frac{1}{2}}-y_{j-\frac{1}{2}}
\end{eqnarray*}

Similarly, for the flux through the southern edge of the control volume $\C_{i,1}~~~~~i=1,...,N$, we have

\begin{equation}
\label{south-appr-fitfin}
\int_{(x_{i-\frac{1}{2}},y_{\frac{1}{2}})}^{(x_{i+\frac{1}{2}},y_{\frac{1}{2}})}\Bigg(m_{21}\frac{\partial \V}{\partial x}+m_{22}\frac{\partial \V}{\partial y}+q\V\Bigg)dx  \approx  \frac{1}{2}y_1\Big[\frac{1}{2}h_i(e+k)-h_i'\Big]\V_{i,1}+\frac{1}{2}h_i'y_1\V_{i+1,1}-\frac{1}{4}y_1h_i(e-k)\V_{i,0}
\end{equation}

where

\begin{eqnarray*}
	e=\frac{1}{2}\sigma_2^2,~~~~~~~~~k=r-\sigma_2^2-\frac{1}{2}\rho\sigma_1\sigma_2  ~~~~~~~h_i'=\frac{1}{2}\rho\sigma_1\sigma_2x_i~~~~~~h_i=x_{i+\frac{1}{2}}-x_{i-\frac{1}{2}}
\end{eqnarray*}


\subsubsection{The fitted L-Multi-Point Flux Approximation method ( with the $1^{st}$ order upwind method ) }


\begin{enumerate}
	\item  Fitted L-MPFA method (with $1^{st}$ order upwind method)
	
	Here the fitted finite volume method is combined with the first order upwind method . Thereby we have:

   \begin{itemize}
		
      \item[] 
      
      For the control volume $\C_{11}$, the western and southern edges of  are (fully) in the degeneracy region. The integrals over the western and southern edges of the control volume $\C_{11}$ are then approximated using the  fitted finite volume (\ref{west-appr-fitfin}) and (\ref{south-appr-fitfin}). The integrals over the eastern and northern edges of the control $\C_{11}$, which are not in the degeneracy region,  are approximated using the L-MPFA method coupled to the upwind method ($1^{st}$ and $2^{nd}$ order).
	
 \begin{eqnarray*}
\label{fitfin-C11}
	\int_{{\mathcal{C}}_{11}}\nabla k(\V) d \C_{11}	&  \approx  & aa_{11}\V_{11}+bb_{11}\V_{21}+cc_{11}\V_{22}+dd_{11}\V_{12}+ee_{11}\V_{01}+\beta\beta_{11}\V_{10}
\end{eqnarray*}
where

\begin{eqnarray*}
	&  & aa_{11}= T_{11}^{22}+T_{34}^{21}+T_{41}^{22}+T_{22}^{12}+l_1\max(f_x^2,0)+h_1\max(f_y^{2},0)-\frac{1}{2}x_1\Big(\frac{1}{2}l_1(a+b)-d_1\Big)\\
	&    & \\
	&    &~~~~~~~~-\frac{1}{2}y_1\Big(\frac{1}{2}h_1(e+k)-h_1'\Big)\\
	&    & \\
	&    & \\
	&    & bb_{11}=T_{33}^{21}+T^{22}_{12}+l_1\min(f_x^2,0)-\frac{1}{2}h_1'y_1~~~~~~~~~~~~~~~~~~cc_{11}=T_{13}^{22}+T_{43}^{22}\\
	&    & \\
	&    & \\
	&    & dd_{11}=T_{23}^{12}+T_{44}^{22}+h_1\min(f_y^{2},0)-\frac{1}{2}d_1x_1~~~~~~~~~~ee_{11}=T_{21}^{12}+\frac{1}{4}l_1x_1(a-b)\\
	&    & \\
	&    & \\
	&    & \beta\beta_{11}=T_{31}^{21}+\frac{1}{4}y_1h_1(e-k)
\end{eqnarray*}

\item[] 

Similarly, for the control volume $\C_{1,j}~~~~j=1,\ldots,N$, only the southern edge is (fully) in the degeneracy  region then the integral over it, is approximated using the fitted finite volume method (\ref{south-appr-fitfin}). The integrals over the eastern, northern and western edges are approximated using the L-MPFA method coupled to the upwind methods($1^{st}$ and $2^{nd}$ order)

\begin{eqnarray}
\label{fitfin-C1j}
	\int_{C_{1,j}} \nabla \mathbf{k}\V d \C_{1 j} & =  & aa_{1,j}\V_{1,j}+bb_{1,j}\V_{2,j}+cc_{1,j}\V_{2,j+1}+dd_{1,j}\V_{1,j+1}+\beta\beta_{1,j}\V_{1,j-1}
	+\alpha\alpha_{1,j}\V_{0,j-1}+ee_{1,j}\V_{0,j} \nonumber\\
\end{eqnarray}

where

\begin{eqnarray*}
	&   & aa_{1,j}=T_{11}^{2,j+1}+T_{34}^{2,j}+T_{41}^{2,j+1}+T_{22}^{1,j+1}-T_{23}^{1,j}-T_{44}^{2,j}+l_j\max(f_x^2,0)+h_1\max(f_y^{j+1},0)-h_1\min(f_y^{j},0)\\
	&    & \\
	&    & ~~~~~~~~~-\frac{1}{2}x_1\Big(\frac{1}{2}l_j(a+b) -d_j\Big)\\
	&    & \\
	&    & \\
	&    & bb_{1,j}=T_{33}^{2,j}+T^{2,j+1}_{12}-T_{43}^{2,j}+l_j\min(f_x^2,0)~~~~~~~~~~~~~~~~~~~~~~~~~~cc_{1,j}=T_{13}^{2,j+1}+T_{43}^{2,j+1}\\
	&    & \\
	&    & \\
	&    & dd_{1,j}=T_{23}^{1,j+1}+T_{44}^{2,j+1} +h_1\min(f_y^{j+1},0)-\frac{1}{2}d_jx_1~~~~~~~~\beta\beta_{1,j}=T_{31}^{2,j}-T_{22}^{1,j}-T_{41}^{2,j}-h_1\max(f_y^{j},0)\\
	&    & \\
	&    & \\
	&    & \alpha\alpha_{1,j}=-T_{21}^{1,j}~~~~~~~~~~~~~~~~~~~~~ee_{1,j}=T_{21}^{1,j+1}+\frac{1}{4}l_jx_1(a-b)
\end{eqnarray*}

\item[] 
Using the same argument as above,     for the control volume $\C_{i,1}~~~i=2,..,N$, the integral over the southern edge is approximated using the fitted finite volume (\ref{south-appr-fitfin}). The integrals over the eastern, northern and western edges are approximated using the L-MPFA method combined with the upwind methods ($1^{st}$ and $2^{nd}$ order)

\begin{eqnarray}
\label{fitfin-Ci1}
	\int_{{\mathcal{C}}_{i,1}}\nabla (\mathbf{k}\V) d \C_{i 1}&   =   & aa_{i,1}\V_{i,1}+bb_{i,1}\V_{i+1,1}+cc_{i,1}\V_{i+1,2}+dd_{i,1}\V_{i,2}+ee_{i,1}\V_{i-1,1}+\alpha\alpha_{i,1}\V_{i-1,0}
	+\beta\beta_{i,1}\V_{i,0} \nonumber\\
\end{eqnarray}

where

\begin{eqnarray*}
	&   & aa_{i,1}=T_{11}^{i+1,2}+T_{34}^{i+1,1}+T_{41}^{i+1,2}+T_{22}^{i,2}-T_{12}^{i,2}-T_{33}^{i,1}+l_1\max(f_x^{i+1},0)
	+h_i\max(f_y^2,0)-l_1\min(f_x^{i},0)\\
	&    & \\
	&    &~~~~~~~~~-\frac{1}{2}y_1\Big(\frac{1}{2}h_i(e+k)-h_i'\Big)\\
	&    & \\
	&    & \\
	&    & bb_{i,1}=T_{33}^{i+1,1}+T^{i+1,2}_{12}+l_1\min(f_x^{i+1},0)-\frac{1}{2}h_i'y_1~~~~~~~~~~~~~~~~~
	cc_{i,1}=T_{13}^{i+1,2}+T_{43}^{i+1,2}\\
	&    & \\
	&    & \\
	&    & dd_{i,1}=T_{23}^{i,2}+T_{44}^{i+1,2}-T_{13}^{i,2}+h_i\min(f_y^2,0)~~~~~~~~~~~~~ee_{i,1}=T_{21}^{i,2}-T_{11}^{i,2}-T_{34}^{i,1}-l_1\max(f_x^{i},0)\\
	&    & \\
	&    & \\
	&    & \alpha\alpha_{i,1}=-T_{31}^{i,1}~~~~~~~~~~~~~~~~\beta\beta_{i,1}=T_{31}^{i+1,1}+\frac{1}{4}y_1h_i(e-k)
\end{eqnarray*}

\end{itemize}

Besides, for the control volume $\C_{ij},~~~~i,j=2,...,N~$, the L-MPFA method is used to approximate the diffusion term and the upwind to approximate the advection term. This leads to the following semi-discrete equation

\begin{equation}
\label{fitlmpfa-up1}
\frac{d\V}{d\tau}=A\V+G(\V)+F
\end{equation} 

where

\begin{eqnarray*}
	\V=\begin{bmatrix}
		\V_{11}\\
		\V_{12}\\
		\vdots\\
		\V_{1N}\\
		\V_{21}\\
		\V_{22}\\
		\vdots\\
		\V_{NN}
	\end{bmatrix}~~~
A=L^{-1}\Bigg(Z+A_L\Bigg)~~~~G(\V)=\beta \Big[\max\Big(\V^*-\V,0\Big)\Big]^{1/k}
\end{eqnarray*}
\\
with F the vector of boundary conditions, $A_L$  a diagonal matrix of size $N\times N$ coming from the discretisation of \eqref{linearterm}. The diagonal elements of $A_L$ are $A_{ii}=h_il_i\lambda$ for $i=1,...,N$ with $\lambda$ given in \eqref{conservation}. The matrix L is also a diagonal matrix of size $N\times N$ whose diagonal elements are $L_{ii}=h_il_i$ for $i=1,\ldots,N$ and

\begin{equation*}
Z=\begin{bmatrix}
D_1 & K_1 & 0_N & \ldots & \ldots & \ldots & \ldots & 0_N\\ 
L_2 & D_2 & K_2 & \ddots &  &  &  & \vdots  \\ 
0_N & L_3 & D_3 & K_3 & \ddots &  & & \vdots\\ 
\vdots & \ddots & L_4 & D_4 & K_4 & \ddots & & \vdots\\ 
\vdots &  & \ddots & \ddots & \ddots & \ddots & \ddots & \vdots \\ 
\vdots &  &  & \ddots& \ddots & \ddots & \ddots & 0_N&  \\ 
\vdots &  &  &  & \ddots & L_{N-1} & D_{N-1} & K_{N-1} \\ 
0_N & \ldots & \ldots & \ldots & \ldots & 0_N & L_N & D_N
\end{bmatrix}   
\end{equation*}

The elements of matrix Z are matrices. $0_N$ is a zeros matrix of size $N$. The matrices $D_i,K_i,L_i$ are tri-diagonal matrices  defined as follows:

\begin{eqnarray*}
	for~~i=1~or~i=N~~~~~~~~~~~~~~~~~~~~~~~~~~~~~~~~~~~~~~~~~~~~~~~~~~~~~~~~~~~~~~~~~~~~
	~~~~~~~~~~~~~&  & \\
	&  & \\
	k=1,\ldots,N~~(D_i)_{kk}=aa_{1,k}~~~~~~~~~k=1,\ldots,N-1~~(D_i)_{k,k+1}=dd_{1,k},~~~~~~~~~~~~k=2,\ldots,N~~(D_i)_{k,k-1}=\beta\beta_{1,k}\\
	&  & \\
	&   & \\
	for~i=1~~~~~~~~~~k=1,\ldots,N~~(K_1)_{kk}=bb_{1,k}~~~~~~~~~~~~~~~~~~~~~~~~~~~~~~k=1,\ldots,N-1~~(K_1)_{k,k+1}=cc_{1,k}\\
	&    &  \\
	&    & \\
for~i=N,~~~~~~~ (L_N)_{11}=ee_{N,1}~~~~k=2,\ldots,N~(L_N)_{kk}=e_{N,k}+\eta_{N,k}~~~~~~~~~~~~k=2,\ldots,N~~~(L_N)_{k,k+1}=\alpha_{N,k}
\end{eqnarray*}

\begin{eqnarray*}
	&  &~~~~~~~~~~~~~~~~~~~~~~~~~~~~~~~~~~~~~~~~~~~~~~~for ~~i=2,\ldots,N-1 \\
	&  & \\
	&  & ~(D_i)_{11}=aa_{i,1}~;~(D_i)_{12}=dd_{i,1};~~~~~(K_i)_{11}=bb_{i,1}~
	;~(K_i)_{12}=cc_{i,1}~~~~(L_i)_{11}=ee_{i,1}\\
	&   & \\
	&   & 
	k=2,\ldots,N~~~~~(D_i)_{kk}=a_{i,k}+\Omega_{i,k};~~~~~~~
	(K_i)_{kk}=b_{i,k}+\Delta_{i,k};~~~~~~~~
	(L_i)_{kk}=e_{i,k}+\eta_{i,k}\\
	&   & \\
	&   & k=2,\ldots,N-1~~~~~(D_i)_{k,k+1}=d_{i,k}+\phi_{i,k};~~~~~~~
	(K_i)_{k,k+1}=c_{i,k}\\
	&   &  \\
	&   & k=2,\ldots,N~~~~~(D_i)_{k,k-1}=\beta_{i,k}+\mu_{i,k};~~~~~~~~~~~~(L_i)_{k,k-1}=\alpha_{i,k}\\
	&    & \\
	&    & \\
	&   & 
\end{eqnarray*}

where  the elements  $aa_{ij},bb_{ij},cc_{ij},dd_{ij},ee_{ij},\beta\beta_{ij}$ are defined in \eqref{fitfin-C11},\eqref{fitfin-C1j},\eqref{fitfin-Ci1}, and the elements $a_{ij},b_{ij},c_{ij},d_{ij},e_{ij},\Omega_{ij}$ $ \Delta_{ij},\beta_{ij},\phi_{ij},\alpha_{ij},\mu_{ij},\eta_{ij}$ are defined in \eqref{mpfa-L} and \eqref{flux-up1}.

\item Fitted Multi-Point Flux Approximation ($2^{nd}$ order upwind)

Similarly, the fitted MPFA-L method deriving from the combination of the L-MPFA method and the $2^{nd}$ order upwind method leads to the following equation :

\begin{equation}
\label{fitlmpfa-up2}
\frac{d\V}{d\tau}=A\V+G(\V)+F
\end{equation} 

where

\begin{eqnarray*}
	\V=\begin{bmatrix}
		\V_{11}\\
		\V_{12}\\
		\vdots\\
		\V_{1N}\\
		\V_{21}\\
		\V_{22}\\
		\vdots\\
		\V_{NN}
	\end{bmatrix}~~~
	A=L^{-1}\Bigg(Y+A_L\Bigg)~~~~G(\V)=\beta \Big[\max\Big(\V^*-\V,0\Big)\Big]^{1/k}
\end{eqnarray*}

with F the vector of boundary conditions. $A_L$ is a diagonal matrix of size $N\times N$ coming from the discretisation of \eqref{linearterm}. The diagonal elements of $A_L$ are $A_{ii}=h_il_i\lambda$ for $i=1,...,N$ with $\lambda$ given in \eqref{conservation}. The matrix L is also a diagonal matrix of size $N\times N$ whose diagonal elements are $L_{ii}=h_il_i$ for $i=1,\ldots,N$

and 

\begin{equation*}
Y=\begin{bmatrix}
H_1 & P_1 &   R_1 & 0_N & \ldots & \ldots & \ldots &  & 0_N & 0_N\\ 
Q_2 & H_2 & P_2& R_2 & 0_N &   &  &   & & 0_N \\
W_3& Q_3 & H_3 & P_3 &  R_3 & 0_N &   &    &  & \vdots  \\
0_N & W_4 & Q_4 & H_4 & P_4 & R_4 & \ddots\\
0_N    &    0_N &  \ddots & \ddots &  \ddots & \ddots & \ddots & \ddots \\
\vdots &   &  \ddots & \ddots & \ddots & \ddots & \ddots & \ddots                                          & \ddots \\
\vdots &   &  & \ddots & \ddots & \ddots & \ddots & \ddots                                          & \ddots & 0_N \\
&   &   &  &  \ddots &  W_{N-2} & Q_{N-2}  & H_{N-2} & P_{N-2}  & R_{N-2} \\
\vdots &  &  &  & & \ddots & W_{N-1} & Q_{N-1} & H_{N-1} & P_{i,N-1}\\
0_N & \dots   &   \ldots  & \ldots  &   & \ldots & 0_N & W_N  & Q_{N}  & H_{N} \\
\end{bmatrix}
\end{equation*}

The elements of matrix Y are matrices. $0_N$ is a zeros matrix of size $N$. The matrices  $H_i~~~~i=1,..,N$ are a penta-diagonal matrices and the matrices $P_i,R_i,Q_i,W_i$ are diagonal matrices.
\\
\\
Furthermore, The $\theta$- Euler method will be applied on the semi-discrete equations (\ref{mpfa-up1}),(\ref{mpfa-up2}), (\ref{fitlmpfa-up1}) and(\ref{fitlmpfa-up2}) for the spatial discretisation.

\section{Time discretization}

Let us consider the ODE stemming from the spatial discretization and given by \eqref{mpfa-up1},\eqref{mpfa-up2},(\ref{fitlmpfa-up1}) and(\ref{fitlmpfa-up2})

\begin{equation*}
\frac{d\V}{d\tau}=A\V+G(\V)+F
\end{equation*}

By using the $\theta$-Euler method for the time discretization, we have

\begin{eqnarray}
\label{euler}
for~m=0,\ldots,M~~~~\frac{\V^{m+1}-\V^m}{\Delta\tau}=\theta\Big(A\V^{m+1}+G(\V^{m+1})+F(t_{m+1}\Big)+(1-\theta)\Big(A\V^m+G(\V^m)+F(t_m\Big) \nonumber\\
\end{eqnarray}
At every time iteration, the nonlinear system  where $\V^{m+1}$   is the solution is solving using the Newton method.
%
%
%
%
%
Note that
\begin{eqnarray*}
	&  & \V^m=\begin{bmatrix}
		\V_{11}(\tau_m)~~
		\V_{12}(\tau_m)~~
		\ldots
		\V_{1N}(\tau_m)~~
		\V_{21}(\tau_m)~~
		\V_{22}(\tau_m)~~
		\ldots
		\V_{2N}(\tau_m)~~
		\ldots
		\V_{N,1}(\tau_m)~~
		\V_{N,2}(\tau_m)~
		\ldots\ldots~
		\V_{NN}(\tau_m)
	\end{bmatrix}^T\\
	&  & \\
	&  & ~~~~~~~~~~~~~~~~~~~~~~~~~~~~~~~~~~~~~~~~~~~~~~~~~~~~~~~~~~~~
	~~~~~~~\tau_m=m\Delta \tau
\end{eqnarray*}
where the time step is $\Delta \tau=\frac{T}{M}$, $T$ being the maturity time.

\end{enumerate}

\section{Numerical experiments}

In this section, we perform some numerical simulations for the L-MPFA method combined to the upwind methods (first and second order) and for the fitted L-MPFA method combined to the upwind methods (first and second order). \\

%
%
%

\subsection{Errors for European call options}

The computational domain of the problem is $\Omega=[0;300]\times [0;300]\times[0;T]$ with T=1/12. The numerical experiments are performed with the strike price $E=100$, the volatilities $\sigma_1=\sigma_2=0.3$, the correlation coefficient $\rho=0.3$ and the risk free interest $r=0.08$\\

Here, by taking $\beta=0$ in (\ref{penpow1}), the L-MPFA method illustrated in the previous sections will be compared to the fitted finite volume method, \cite{wang2004novel}, and the   fitted O-MPFA  methods for pricing  multi-asset options for pricing options introduced in \cite{rock2019mpfa}. The relative error will be computed with respect to the analytical solution  of the Black-Scholes PDE  defined in \cite{haug2007complete} as follows
\begin{eqnarray}
\label{analsol}
C(S_1,S_2,K,T) & = &  S_1e^{-rT}M(y_1,d;\rho_1)+S_2e^{-rT}M(y_2,-d+\sigma\sqrt{T};\rho_2) \nonumber\\
&   &   \\
&    & -Ke^{-rT}\times\left(1-M(-y_1+\sigma_1\sqrt{T},-y_2+\sigma_2\sqrt{T};\rho)\right) \nonumber
\end{eqnarray} 

where

\begin{eqnarray*}
	&  & d=\frac{ln(S_1/S_2)+(b_1-b_2+\sigma_1^2/2)T}{\sigma\sqrt{T}},\\
	&   & \\
	&   & y_1  =  \frac{ln(S_1/K)+(b_1+\sigma_1^2/2)T}{\sigma_1\sqrt{T}},~~~~~~y_2=\frac{ln(S_2/K)+(b_1+\sigma_2^2/2)T}{\sigma_2\sqrt{T}},\\
	&   &   \\
	&   & \sigma=\sqrt{\sigma_1^2+\sigma_2^2-2\rho\sigma_1\sigma_2},~~~~~\rho_1=\frac{\sigma_1-\rho\sigma_2}{\sigma}~~~~~~~\rho_2=\frac{\sigma_2-\rho\sigma_1}{\sigma}
\end{eqnarray*}

and

\begin{equation*}
M(a,b;\rho)=\frac{1}{2\pi\sqrt{1-\rho^2}}\int_{-\infty}^a\int_{-\infty}^bexp\left(-\frac{x^2-2\rho xy+y^2}{2(1-\rho^2)}\right)dxdy.
\end{equation*}

\begin{figure}[hbtp]
	\centering
    \includegraphics[scale=0.5]{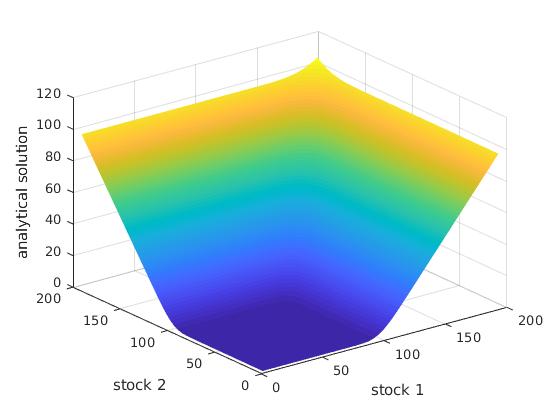}
	\caption{Analytical solution}
	\label{anasol}
\end{figure}

The solution using the L-MPFA coupled to the $2^{nd}$ order upwind method is illustrated as below

\begin{figure}[hbtp]
	\centering
	\includegraphics[scale=0.5]{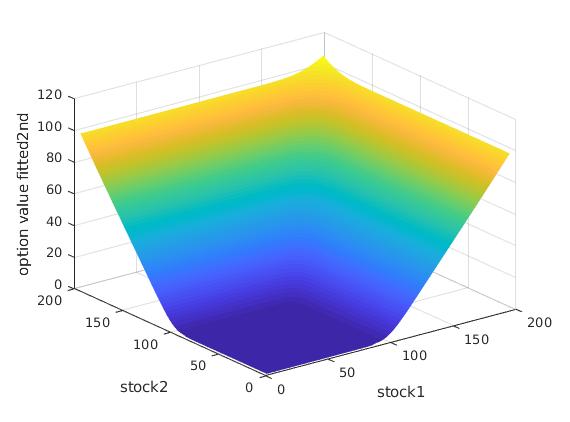}
	\caption{L-MPFA -upwind $2^{nd}$ order}
	\label{anasol}
\end{figure}

 The  $L^2$-norm is used to compute the error is  

\begin{equation}
\label{error-l2}
err=\frac{\sqrt{\sum_{i,j=1}^N meas(\C_{ij}) \big(\V_{ij}-V_{ij}^{ana}\big)^2}}{\sqrt{\sum_{i=1}^n meas(\C_{ij}) \big(V_{ij}^{ana}\big)^2}}
\end{equation}

where $\V$ is the numerical solution, $V^{ana}$ the analytical solution and $meas(\C_{ij})$ is the measure of the control volume $\C_{ij}$.
This gives the following tables:

\begin{table}[!h]
	\centering
	\begin{tabular}{|c||c|c|c|c|c|}
		\hline
		& Fitted fin vol & O-MPFA-$1^{st}$
		upw
		& O-MPFA-$2^{nd}$ upw &  fit O-MPFA-$1^{st}$ upw & fit  O-MPFA -$2^{nd}$ upw\\
		\hline 
		$50\times 50$ & 0.0317 &   0.0224 & 0.0225 &  0.0212  & 0.0212 \\
		\hline 
		$70\times 70$ & 0.0329 &  0.0248 & 0.0248  &  0.0238 & 0.0238 \\
		\hline
		$85\times 85$ & 0.0327 & 0.0260     &    0.0260   & 0.0251    &  0.0251\\
		\hline
	\end{tabular}
	\caption{Table of errors}
	\label{table:euro1}
\end{table}

\begin{table}[!h]
	\centering
	\begin{tabular}{|c||c|c|c|c|c|}
		\hline
		& L-MPFA-$1^{st}$upw 
		& L-MPFA-$2^{nd}$ upw &  fit L-MPFA-$1^{st}$ upw & fit L-MPFA -$2^{nd}$ upw\\
		\hline 
		$50\times 50$		 & 0.0048 &   0.0049 & 0.0048 &  0.0047  \\
		\hline 
		$70\times 70$		 & 0.0041 &  0.0041 & 0.0041  &  0.0041  \\
		\hline
		$85\times 85$		& 0.0040 & 0.0040     &    0.0040   & 0.0040  \\
		\hline
	\end{tabular}
	\caption{Table of errors}
	\label{table:euro2}
\end{table}

As we can observe in  Table \ref{table:euro1} and Table \ref{table:euro2},  the new fitted L-MPFA method is more accurate than the  fitted O-MPFA  method developed in \cite{rock2019mpfa}  and the standard fitted finite volume method developed in \cite{huang2006fitted}.

\subsection{Errors for American put options}

Since there is no analytical solution for the power penalty problem (\ref{penpow1}) for pricing American put options, and the numerical solution given by the fitted L-MPFA coupled to $2^{nd}$ order upwind method is more accurate when  pricing European options (see Table \ref{table:euro1} and Table \ref{table:euro2}),  we have  chosen for reference solution  or ''exact solution'' the numerical solution given by the fitted L-MPFA coupled to $2^{nd}$ order upwind method  with $dt=T/256$. The  relative error   of all the numerical methods used in this this study will be performed with respect to this reference solution.

\begin{figure}[hbtp]
	\centering
	\includegraphics[scale=0.5]{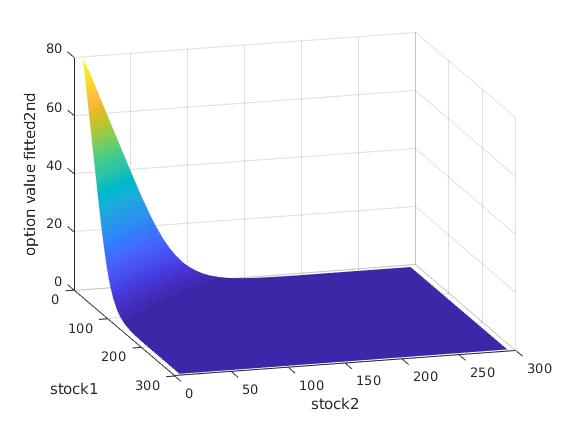}
	\caption{Reference solution}
	\label{solref}
\end{figure}

%
%

For the  numerical simulations below, the computational domain of the problem is $\Omega = [0; 300] \times [0; 300] \times  [0; T ]$ with
$T = 1/6, K = 100$, the volatilities $\sigma_1 = \sigma_2 = 0.3$. The correlation coefficient is $\rho = 0.3$ , the risk free interest
r = 0.08. The penalty parameter $\beta=256$ and the power penalty $k=1/2$.

\begin{table}[!h]
	\centering
	\begin{tabular}{|c||c|c|c|c|c|c|}
		\hline
		& Fitted Fin vol & L-MPFA-$1^{st}$upw 
		& L-MPFA-$2^{nd}$ upw &  fit L-MPFA-$1^{st}$ upw & fit L-MPFA -$2^{nd}$ upw\\
		\hline 
		$50\times 50$ &	0.0616	 & 0.0610 & 0.0583 & 0.0611 &  0.0584 \\
		\hline 
		$60\times 60$ &	0.0277	 & 0.0278 & 0.0276  & 0.0278 & 0.0277 \\
		\hline
		$70\times 70$	& 0.0184	& 0.0183 &  0.0182  &  0.0182     & 0.0180 \\
		\hline
		$ 80\times 80 $  & 0.0104  & 0.0100  & 0.0098  &  0.0097 & 0.0095 \\
		\hline
	\end{tabular}
	\caption{Table of errors for $\Delta \tau = T/64$}
	\label{table:amer}
\end{table}

\begin{table}[!h]
	\centering
	\begin{tabular}{|c||c|c|c|c|c|c|}
		\hline
		& Fitted Fin vol & L-MPFA-$1^{st}$upw 
		& L-MPFA-$2^{nd}$ upw &  fit L-MPFA-$1^{st}$ upw & fit L-MPFA -$2^{nd}$ upw\\
		\hline 
		$50\times 50$	 & 0.0599 & 0.0520  & 0.0476 & 0.0522 & 0.0459 \\
		\hline 
		$60\times 60$		 & 0.0227 & 0.0265   & 0.0249  & 0.0241 & 0.0220 \\
		\hline
		$70\times 70$		&  0.0136  &  0.0148    & 0.0146      &  0.0146 & 0.0144 \\
		\hline
		$ 80\times 80 $   & 0.0087  & 0.0068  &  0.0065 & 0.0062  & 0.0059 \\
		\hline
	\end{tabular}
	\caption{Table of errors for $\Delta \tau = T/128$}
	\label{table:amer1}
\end{table}

 Again  we can observe  in Table \ref{table:amer}  and Table \ref{table:amer1},   the novel fitted L-MPFA coupled to the $2^{nd}$ order upwind method  is more accurate than  the  standard fitted finite volume by \cite{huang2006fitted}.

\section{Conclusion}

In this paper, the L-MPFA methods  have been introduced to approximate the diffusion term of the Black-Scholes PDE. The upwind methods ($1^{st}$ and $2^{nd}$ order) are used for  space discretization of  the convection term appearing in the two dimensional  Black-Scholes PDE. We have provided a novel scheme called the fitted L-MPFA method to handle the degeneracy of the Black Scholes PDE by combining   the fitted finite volume and the L-MPFA method coupled to the upwind methods.  Numerical experiments are performed and comparison between the L-MPFA methods, the O-MPFA methods by \cite{rock2019mpfa} and the fitted finite methods by \cite{huang2006fitted} are performed. The results have shown that the fitted L-MPFA method coupled to the  $2^{nd}$ order upwind method is more accurate than the fitted finite volume by \cite{huang2006fitted} and the O-MPFA  method by \cite{rock2019mpfa} for pricing Europeans and American options.

\section*{Acknowledgement}

This work was supported by the Robert Bosch Stiftung through the AIMS ARETE Chair programme (Grant No 11.5.8040.0033.0).
%
	
%
\end{document}